\documentclass[a4paper,11pt]{article}
\usepackage[utf8]{inputenc}
\usepackage{hyperref}
\usepackage{amsfonts}
\usepackage{amssymb}
\usepackage{amsthm}
\usepackage{epsf, epsfig}
\usepackage{amssymb,latexsym}
\usepackage{url,tikz,yhmath,appendix,booktabs}
\usepackage{color}
\usepackage[normalem]{ulem}
\usepackage{cancel} %\cancel{}

\usepackage{subfigure}
\usepackage{graphicx}
\usepackage{lipsum}
\usepackage{pinlabel}
\usepackage{authblk}%%%for authors

\newcommand{\red}{\textcolor{red}}

\newtheorem{theo}{Theorem}%[section]
\newtheorem{pro}{Proposition}
\newtheorem{cor}{Corollary}

\newtheorem{prob}{Problem}

\usetikzlibrary{backgrounds}
\usetikzlibrary{topaths,calc}
\usetikzlibrary{backgrounds}
\usetikzlibrary{arrows, patterns}
\usetikzlibrary{arrows.meta}
\usetikzlibrary{patterns.meta}
\usetikzlibrary{shapes,shapes.geometric,shapes.misc}
\pgfdeclarelayer{edgelayer}
\pgfdeclarelayer{nodelayer}
\pgfsetlayers{background,edgelayer,nodelayer,main}
\tikzstyle{none}=[inner sep=0mm]
\tikzstyle{every loop}=[]

\tikzstyle{dotted}=[dash pattern=on \pgflinewidth off 2pt]
\tikzstyle{dashed}=[dash pattern=on 3pt off 3pt]

\newcommand \tikzp[2]
{
	\begin{center}
		\begin{tikzpicture}[scale=#1]
			#2
		\end{tikzpicture}
	\end{center}
}

\tikzstyle{new style 0}=[fill=black, draw=black, shape=circle]
\tikzstyle{red style 1}=[fill=red, draw=black, shape=circle]
\tikzstyle{blue style 2}=[fill=blue, draw=black, shape=circle]
\tikzstyle{white style 4}=[fill=white, draw=black, shape=circle]
\tikzstyle{bklack style 5}=[fill=black, draw=black, shape=rectangle]
\tikzstyle{red style 3}=[fill=red, draw=black, shape=rectangle]
\tikzstyle{yellow style 7}=[fill=yellow, draw=black, shape=rectangle]
\tikzstyle{new style 8}=[fill={rgb,255: red,0; green,132; blue,0}, draw={rgb,255: red,0; green,131; blue,0}, shape=circle]

\tikzstyle{new edge style 0}=[-]
\tikzstyle{new edge style 1}=[-, draw=red]
\tikzstyle{new edge style 2}=[-, draw=blue]
\tikzstyle{new edge style 3}=[-, draw={rgb,255: red,0; green,156; blue,0}]

\tikzstyle{cblue}=[circle, draw, thin,fill=blue!20, scale=0.5]
\newcommand\idf[2]
{
_{#1\bullet  #2}
}

\newcommand{\proofend}{{\hfill$\Box$} \vspace{0.3cm}}

\def\rebibitem {\bibitem}  %final version

\begin{document}

\title{A study on $T$-equivalent graphs
}

\author[1]{\small Fengming Dong\thanks{Email: fengming.dong@nie.edu.sg;
		donggraph@163.com}}
\author[2]{\small Meiqiao Zhang\thanks{
		Corresponding Author. Email: %nie21.zm@e.ntu.edu.sg and 
		meiqiaozhang@xmu.edu.cn;
		meiqiaozhang95@163.com}}

\affil[1]{\footnotesize 
	National Institute of Education,
	Nanyang Technological University, 
	Singapore}

\affil[2]{\footnotesize School of Mathematical Sciences, Xiamen  University, China}

\date{}

\maketitle

\begin{abstract}
In his article [{\it J. Comb. Theory Ser. B} {\bf 16} (1974), 168--174],
Tutte called two graphs 
$T$-equivalent 
(i.e., codichromatic) if they have the same Tutte polynomial and  
showed that graphs $G$ and $G'$ are
$T$-equivalent 
if $G'$ is obtained from $G$ 
by flipping a rotor
(i.e., replacing it by its mirror) of order at most $5$, 
where a rotor of order $k$ in $G$ is an induced subgraph $R$ having 
an automorphism $\psi$ with a vertex orbit $\{\psi^i(u): i\ge 0\}$ of size $k$ 
such that every vertex of $R$ 
is only adjacent to vertices in $R$
unless it is in this vertex orbit.
In this article, we first show the above result due to 
Tutte can be extended to a rotor $R$ of order 
$k\ge 6$ 
if the subgraph of $G$ induced by all those edges of $G$ 
which are not in $R$ satisfies certain conditions.  
Also, we provide  a new method for 
 generating infinitely many non-isomorphic $T$-equivalent pairs of graphs.
\iffalse 
$\chi$-equivalent graphs (i.e., graphs 
having the same chromatic polynomial)
by studying $\chi$-equivalent graphs $R$ and $R'$
and edges $e$ in $R$ and $e'$ in $R'$ 
such that $R\backslash e$ and $R'\backslash e'$ are also cochromatic, where $R\backslash e$ is the graph obtained from $R$ 
by removing $e$.\fi 
\end{abstract}



%\Large  \tableofcontents 

\section{Introduction}\label{intro}

The graphs considered in this article are multigraphs
which may have loops and parallel edges
unless it is mentioned. 
For any graph $G$, 
let $V(G)$ (or simply $V$)
and $E(G)$ (or simply $E$) be its vertex set and edge set.
The Tutte polynomial of $G$ is defined as follows (see~\cite{tutte94}):
\begin{equation}\label{tuttep}
T_G(x,y)=\sum_{A\subseteq E}(x-1)^{r(E)-r(A)}(y-1)^{|A|-r(A)},
\end{equation}
where $r(A)=|V|-c(A)$ and $c(A)$ is the number of components 
of the spanning subgraph $(V,A)$ of $G$.
The Tutte polynomial is an invariant so fundamental and all-encompassing that many well-known graph polynomials and research objects originated from other disciplines have been shown to be its specializations~\cite{Dong,sokal,tutte2,welsh}. Over the years, the Tutte polynomial has been extensively studied from various perspectives, especially on the combinatorial explanations of its specific coefficients and evaluations~\cite{bollobas,tutte94}. 
All of this reveals that the Tutte polynomial contains far more information  than previously understood, highlighting its significant research value.

\def \tsim {\stackrel T\sim}

In this paper, we focus on the study of {\it $T$-equivalent graphs}, or {\it codichromatic graphs} by Tutte~\cite{tutte}, i.e., two graphs $G_1$ and $G_2$ having the same Tutte polynomial, written as 
$G_1\tsim G_2$. Clearly, two isomorphic graphs 
are $T$-equivalent. Also, Brylawski and Oxley~\cite{oxley1} showed that two graphs having isomorphic cyclic matroids are $T$-equivalent. For further developments, we refer the readers to~\cite{boll,Bonin,bry}. In particular, some operations have been found useful in constructing $T$-equivalent pairs. A well-known example is the \textit{Whitney twist}~\cite{whit},
which changes a graph to another one
by flipping a subgraph 
 at a vertex-cut of size $2$. See such a pair $G,G'$ in Figure~\ref{f1-0}, where $G'$ is obtained from $G$ by a Whitney twist with respect to the cut-set $\{u_1,u_2\}$ of $G$.

\iffalse
Two graphs $G_1$ and $G_2$ having the same Tutte polynomial
are called {\it codichromatic graphs} by Tutte~\cite{tutte}, i.e., 
{\it $T$-equivalent graphs}, 
written as $G_1\tsim G_2$.
It is trivial that two isomorphic graphs 
are $T$-equivalent. 
If two graphs have isomorphic cyclic matroids, then 
they are also $T$-equivalent (see~\cite{oxley1}).
A well-known operation for constructing 
such a pair of graphs 
is the Whitney twist~\cite{whit} 
which changes a graph to another one
by flipping a subgraph 
 at a vertex-cut of size $2$.
 %\footnote{} 
 An example for such two graphs $G$ and $G'$
 is shown in Figure~\ref{f1-0}, where $G'$ is obtained from $G$ by a Whitney twist with respect to the cut-set $\{u_1,u_2\}$.
 \fi

\tikzstyle{cblue}=[circle, draw, thin,fill=blue!20, scale=0.5]
\begin{figure}[!ht]
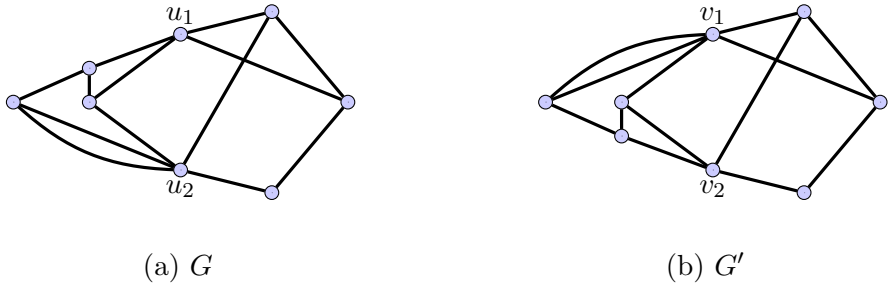

	\tikzp{1}
	{
		\foreach \place/\y in {{(0,0.9)/1}, {(1.2,1.2)/2},{(2.2,0)/3}, {(1.2,-1.2)/4}, {(0,-0.9)/5}, {(-1.2,0)/6},{(-2.2,0)/7},{(-1.2,0.45)/8}}   
		\node[cblue] (b\y) at \place {};
		
		\filldraw[black] (b1) circle (0pt)node[anchor=south] {$u_1$};
		\filldraw[black] (b2) circle (0pt)node[anchor=south] {};
		\filldraw[black] (b3) circle (0pt)node[anchor=south] {};
		\filldraw[black] (b4) circle (0pt)node[anchor=south] {};
		\filldraw[black] (b5) circle (0pt)node[anchor=north] {$u_2$};
		\filldraw[black] (b6) circle (0pt)node[anchor=east] {};
		\filldraw[black] (b7) circle (0pt)node[anchor=north] {};
		\filldraw[black] (b8) circle (0pt)node[anchor=west] {};		
		
		\draw[black, very thick] (b1) -- (b2) -- (b3)--(b4) -- (b5) -- (b7)--(b8)--(b6) --(b1) -- (b8); 
		\draw[black, very thick] (b6) -- (b5) -- (b2);		 
		\draw[black, very thick] (b1) -- (b3);
		\draw [black, very thick, bend left=20] (b5) to (b7);

		\foreach \place/\y in {{(0+7,0.9)/1}, {(1.2+7,1.2)/2},{(2.2+7,0)/3}, {(1.2+7,-1.2)/4}, {(0+7,-0.9)/5}, {(-1.2+7,0)/6},{(-2.2+7,0)/7},{(-1.2+7,-0.45)/8}}   
		\node[cblue] (c\y) at \place {};
		
		\filldraw[black] (c1) circle (0pt)node[anchor=south] {$v_1$};
		\filldraw[black] (c2) circle (0pt)node[anchor=south] {};
		\filldraw[black] (c3) circle (0pt)node[anchor=south] {};
		\filldraw[black] (c4) circle (0pt)node[anchor=south] {};
		\filldraw[black] (c5) circle (0pt)node[anchor=north] {$v_2$};
		\filldraw[black] (c6) circle (0pt)node[anchor=east] {};
		\filldraw[black] (c7) circle (0pt)node[anchor=north] {};
		\filldraw[black] (c8) circle (0pt)node[anchor=west] {};		
		
		\draw[black, very thick] (c1) -- (c2) -- (c3)--(c4) -- (c5) -- (c6)--(c8)--(c7) --(c1) -- (c6); 
		\draw[black, very thick] (c8) -- (c5) -- (c2);
		\draw[black, very thick] (c1) -- (c3);
		\draw [black, very thick, bend right=20] (c1) to (c7);
			}
		
	{}\hfill \hspace{-0.2cm} (a) $G$ \hspace{5.7 cm} (b) $G'$\hfill {}	
			
	\caption{$G'$ is obtained from $G$ by a Whitney twist}
	\label{f1-0}
\end{figure}
 
%Tutte~\cite{tutte} constructed nontrivial $T$-equivalent graphs with connectivity up to $5$. 
%Brylawski~\cite{bry}  and Bollob\'as, Pebody and Riordan~\cite{boll}generalized Tutte's construction to arbitrarily high connectivity.
%of $T$-equivalent graphs.
Another operation is flipping rotors due to Tutte~\cite{tutte}.
Assume that $R$ is a graph 
with an automorphism $\psi$.
For any vertex $x$ in $R$, the set  
$\{\psi^i(x): i\ge 0\}$ is called a {\it vertex orbit}
of $\psi$.
% and $x$ is called a {\it fixed vertex} of $\psi$ if $\psi(x)=x$. 
If $R$ is a subgraph of a graph $G$,
a subset $B$ of $V(R)$
is called a {\it border} of $R$ in $G$ if
every edge in $G$ incident with some vertex  in 
$V(R)\setminus B$ must be an edge in $R$.
We call $R$ a {\it rotor} of $G$ with a border 
$B$ if $B$ is a vertex orbit of 
some automorphism $\psi$.
Tutte~\cite{tutte} showed that if %$k\le 5$ and 
$G$ is a graph containing a rotor $R$   
with a border $B$ of size at most $5$,  
then  $G\tsim G'$,
where $G'$ is the graph 
obtained from $G$ by flipping $R$ along 
its border $B$, i.e., 
by replacing $R$ by its mirror image. 
A formal expression is given below.
%We will express Tutte's result below. 

\iffalse 
Let $G$ and $W$ be vertex-disjoint graphs. 
For distinct vertices $u_1,u_2\in V (G)$
and distinct vertices $w_1,w_2\in V (W)$, 
let $G(u_1, u_2) \sqcup W(w_1, w_2)$
be the graph obtained from $G$ and $W$ by 
identifying $u_1$ and $w_1$ as a vertex 
and identifying $u_2$ and $w_2$ as a vertex. 
We say that the graph $G(u_2, u_1) \sqcup W(w_1, w_2)$
is obtained 
from the graph $G(u_1, u_2) \sqcup W(w_1, w_2)$
by a Whitney twist \cite{whit}.
As these two graphs 
%The two graphs $G(u_2, u_1) \sqcup W(w_1, w_2)$
%and  $G(u_1, u_2) \sqcup W(w_1, w_2)$
have the same cyclic matroid \cite{whit},
they are $T$-equivalent % have the same Tutte Polynomial 
\cite{oxley1}.

\begin{theo}
\label{whitney twist}
The two graphs $G(u_1, u_2) \sqcup W(w_1, w_2)$ 
and 
$G(u_2, u_1) \sqcup W(w_1, w_2)$
have the same Tutte polynomial.
\end{theo}

Theorem~\ref{whitney twist} motivates 
the study of the following problem:
\fi

For any positive integer $k$,
let $[k]$ denote the set $\{1,2,\dots,k\}$.
Given any vertex-disjoint graphs $G$ and 
$W$ with $\{u_i: i\in [k]\}\subseteq 
V(G)$ 
and $\{w_i: i\in [k]\} \subseteq V(W)$,
let $G(u_1,u_2,\dots,u_k)\sqcup 
W(w_1,w_2,\dots,w_k)$ denote the graph obtained 
from $G$ and $W$ by identifying 
$u_i$ and $w_i$ as a new vertex 
for all $i\in [k]$.
An example of $G(u_1,u_2,u_3)\sqcup
W(w_1,w_2,w_3)$ is shown in Figure \ref{f1}.

\tikzstyle{cblue}=[circle, draw, thin,fill=blue!20, scale=0.5]
\begin{figure}[!ht]
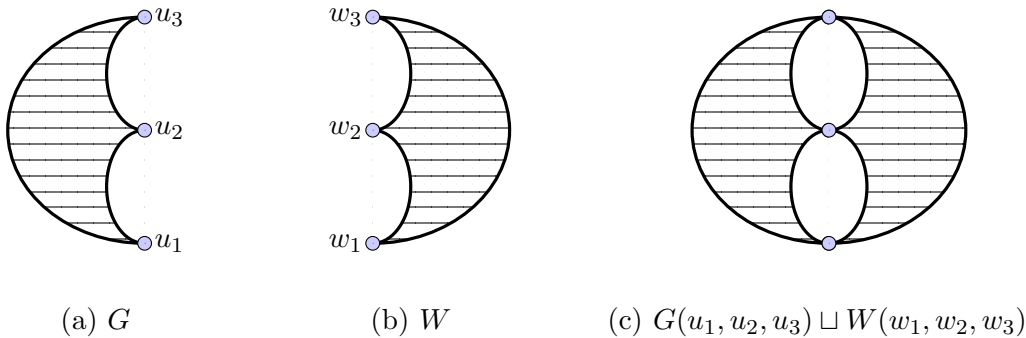

	\tikzp{1}
	{
		\draw[dashed, pattern={Lines[angle=45,distance=6pt]}](0,1.5) 	arc(90:270:1.8cm and 1.5cm);
		\draw[very thick](0,1.5) 	arc(90:270:1.8cm and 1.5cm);		
		\draw[fill=white, very thick](0,1.5) 	arc(90:270:0.5cm and 0.75cm);
		\draw[fill=white, very thick](0,0) 	arc(90:270:0.5cm and 0.75cm);
						
		\foreach \place/\y in {{(0,0)/1}, {(0,1.5)/2},{(0,-1.5)/3}}  
		\node[cblue] (b\y) at \place {};
		
		\filldraw[black] (b1) circle (0pt)node[anchor=west] {$u_2$};
		\filldraw[black] (b2) circle (0pt)node[anchor=west] {$u_3$};
		\filldraw[black] (b3) circle (0pt)node[anchor=west] {$u_1$};

		\draw[dashed, pattern={Lines[angle=45,distance=6pt]}](3,-1.5) 	arc(-90:90:1.8cm and 1.5cm);
\draw[very thick](3,-1.5) 	arc(-90:90:1.8cm and 1.5cm);		
\draw[fill=white, very thick](3,-1.5) 	arc(-90:90:0.5cm and 0.75cm);
	\draw[fill=white, very thick](3,0) 	arc(-90:90:0.5cm and 0.75cm);
						
		\foreach \place/\y in {{(3,0)/1}, {(3,1.5)/2},{(3,-1.5)/3}}  
		\node[cblue] (c\y) at \place {};
		
		\filldraw[black] (c1) circle (0pt)node[anchor=east] {$w_2$};
		\filldraw[black] (c2) circle (0pt)node[anchor=east] {$w_3$};
		\filldraw[black] (c3) circle (0pt)node[anchor=east] {$w_1$};

		\draw[dashed,pattern={Lines[angle=45,distance=6pt]}](9,1.5) 	arc(90:270:1.8cm and 1.5cm);
		\draw[very thick](9,1.5) 	arc(90:270:1.8cm and 1.5cm);		
		\draw[fill=white, very thick](9,1.5) 	arc(90:270:0.5cm and 0.75cm);
		\draw[fill=white, very thick](9,0) 	arc(90:270:0.5cm and 0.75cm);
						
		\foreach \place/\y in {{(9,0)/1}, {(9,1.5)/2},{(9,-1.5)/3}}  
		\node[cblue] (b\y) at \place {};
		
		\filldraw[black] (b1) circle (0pt)node[anchor=west] {};
		\filldraw[black] (b2) circle (0pt)node[anchor=west] {};
		\filldraw[black] (b3) circle (0pt)node[anchor=west] {};

		\draw[dashed,pattern={Lines[angle=45,distance=6pt]}](9,-1.5) 	arc(-90:90:1.8cm and 1.5cm);
\draw[very thick](9,-1.5) 	arc(-90:90:1.8cm and 1.5cm);		
\draw[fill=white, very thick](9,-1.5) 	arc(-90:90:0.5cm and 0.75cm);
	\draw[fill=white, very thick](9,0) 	arc(-90:90:0.5cm and 0.75cm);
						
		\foreach \place/\y in {{(9,0)/1}, {(9,1.5)/2},{(9,-1.5)/3}}  
		\node[cblue] (c\y) at \place {};
		
		\filldraw[black] (c1) circle (0pt)node[anchor=east] {};
		\filldraw[black] (c2) circle (0pt)node[anchor=east] {};
		\filldraw[black] (c3) circle (0pt)node[anchor=east] {};		
}
	{}\hfill \hspace{1.5cm} (a) $G$ \hspace{2.9 cm} (b) $W$  \hspace{1.8 cm} (c) $G(u_1,u_2,u_3)\sqcup
	W(w_1,w_2,w_3)$ \hfill {}			
  \caption{Graphs $G,W$ and $G(u_1,u_2,u_3)\sqcup
W(w_1,w_2,w_3)$}
\label{f1}
\end{figure}

\begin{theo}[Tutte~\cite{tutte}]\label{theo-tutte}
Let $R$ and $W$ be connected graphs.
If $w_1,\dots,w_k$ are distinct vertices in $W$, and $\{u_i: i\in [k]\}$
is a vertex orbit of an automorphism $\psi$ in $R$ such that
$\psi(u_k)=u_{1}$ 
and $\psi(u_i)=u_{i+1}$
for all $i\in [k-1]$, 
where $k\in [5]$, then  
$R(u_{1},\dots,u_k)\sqcup W(w_1,\dots,w_k)$
and $R(u_{k},\dots,u_1)\sqcup W(w_1,\dots,w_k)$ are 
$T$-equivalent.
\end{theo}

Tutte~\cite{tutte} also raised the question of whether Theorem~\ref{theo-tutte} can be extended to the case of $k\ge6$. Although Foldes~\cite{folders}
gave a negative answer to this question by constructing a counterexample pair $(R,W)$ for every $k\ge 6$, our main result in Section~\ref{secmain} (see Theorem~\ref{more-rotors})  provides a sufficient condition on $W$ such that the extension works for any $R$ with $k\ge 6$ that $W$ contains a subgraph with two automorphisms and each of them has a vertex orbit satisfying certain conditions.
 
In order to prove Theorem~\ref{more-rotors},
we first prepare some preliminary results in Section~\ref{sec1} and present in Section~\ref{sec2}, a necessary 
and sufficient condition on graphs $G$ and $H$ 
with $\{u_i: i\in [k]\}\subseteq V(G)$ 
and $\{v_i: i\in [k]\}\subseteq V(H)$
such that  the two graphs 
$G(u_1,u_2,\dots,u_k)\sqcup 
W(w_1,w_2,\dots,w_k)$ and $ H(v_1,v_2,\dots,v_k)\sqcup 
W(w_1,w_2,\dots,w_k)$ are $T$-equivalent 
for an arbitrary graph $W$, 
where
$w_1,\dots,w_k$ are distinct vertices in $W$.

In addition to the graph operations
mentioned above, we are also going to develop an operation to construct $T$-equivalent pairs through the recursive expression of the Tutte polynomial (see (\ref{tuttep2})).
For any graph $G=(V,E)$ and $e=uv\in E$, let  $G\backslash e$ be the graph obtained from $G$ 
by deleting $e$, and let $G/e$ be the graph obtained from $G$ by \textit{contracting} $e$, i.e., the graph obtained from 
$G\backslash e$ by identifying $u$ and $v$ as one vertex.
It is known that $T_G(x,y)$ can also be determined 
by the following recursive expression~\cite{oxley1}:
\begin{equation}\label{tuttep2}
T_G(x,y)=
\left \{ 
\begin{array}{ll}
1, &\mbox{if }E=\emptyset;\\
xT_{G/e}(x,y), &\mbox{if }e \mbox{ is a bridge}; \\
yT_{G\backslash e}(x,y), &\mbox{if }e \mbox{ is a loop}; \\
T_{G/e}(x,y)+T_{G\backslash e}(x,y)
, &\mbox{otherwise}. \\
\end{array}
\right.
\end{equation}

Let $\Phi$ be the set of quaternions $(G,e,H,f)$, 
where $G$ and $H$ are loopless graphs
and $e$ and $f$ are edges in $G$ and $H$ respectively 
such that 
$G\backslash e\tsim H\backslash f$ 
and $G/e\tsim H/f$ hold.
Clearly, for any $(G,e,H,f)\in \Phi$, 
$G\tsim H$ holds.
We will consider the following problem.

\begin{prob}\label{prob2}
Searching for members $(G,e,H,f)$ of $\Phi$.
\end{prob}

\iffalse
In this section, we present a method 
for generating members in $\Phi$
from a known member in $\Phi$.
\fi
\iffalse
Two loopless graphs $G$ and $H$ are said to be 
{\it isomorphic}, written as $G\cong H$,
if there is 
a bijective mapping $\eta$ from $V(G)$ to $V(H)$ 
such that 
$\epsilon_G(u,v)=\epsilon_H(\eta(u),\eta(v))$ holds
for all pairs of vertices $u,v\in V(G)$,
where $\epsilon_G(u,v)$ is the number of edges in $G$ 
joining $u$ and $v$.
\fi

Observe that $\Phi$ contains a subfamily $\Phi'$ 
which consists of those quaternions $(G,e,H,f)$, 
where $G$ and $H$ are connected and loopless graphs
and $e$ and $f$ are edges in $G$ and $H$ respectively
 such that $G\backslash e\cong H\backslash f$
 and $G/e\cong H/f$.

\tikzstyle{cblue}=[circle, draw, thin,fill=blue!20, scale=0.5]
\begin{figure}[!ht]
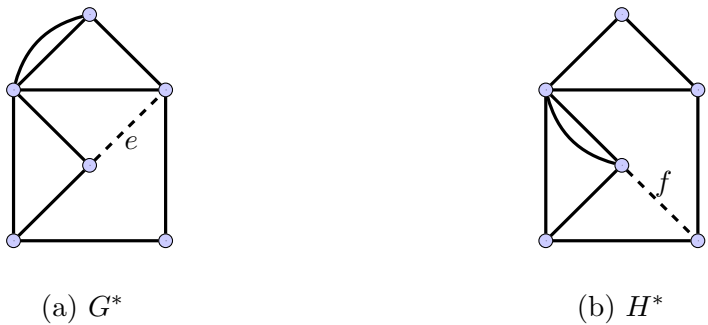

	\tikzp{1}
	{
		\foreach \place/\y in {{(0,0)/1}, {(1,1)/2},{(1,-1)/3}, {(-1,-1)/4}, {(-1,1)/5}, {(0,2)/6}}   
		\node[cblue] (b\y) at \place {};
		
		\filldraw[black] (b1) circle (0pt)node[anchor=south] {};
		\filldraw[black] (b2) circle (0pt)node[anchor=south] {};
		\filldraw[black] (b3) circle (0pt)node[anchor=south] {};
		\filldraw[black] (b4) circle (0pt)node[anchor=south] {};
		\filldraw[black] (b5) circle (0pt)node[anchor=north] {};
		\filldraw[black] (b6) circle (0pt)node[anchor=east] {};

				\node [style=none] (cap1) at (0.55,0.3)  {$e$};
				
		\draw[black, very thick] (b2) -- (b3) -- (b4) -- (b1) --(b5) --(b6)--(b2)--(b5)--(b4);		 
		\draw[black, very thick, dashed] (b1) -- (b2);
		\draw [black, very thick, bend left=30] (b5) to (b6);

		\foreach \place/\y in {{(0+7,0)/1}, {(1+7,1)/2},{(1+7,-1)/3}, {(-1+7,-1)/4}, {(-1+7,1)/5}, {(0+7,2)/6}}   
		\node[cblue] (c\y) at \place {};
		
		\filldraw[black] (c1) circle (0pt)node[anchor=south] {};
		\filldraw[black] (c2) circle (0pt)node[anchor=south] {};
		\filldraw[black] (c3) circle (0pt)node[anchor=south] {};
		\filldraw[black] (c4) circle (0pt)node[anchor=south] {};
		\filldraw[black] (c5) circle (0pt)node[anchor=north] {};
		\filldraw[black] (c6) circle (0pt)node[anchor=east] {};

				\node [style=none] (cap1) at (0.55+7,-0.25)  {$f$};
				
		\draw[black, very thick] (c4) -- (c3) -- (c2) -- (c6) --(c5) --(c4)--(c1)--(c5)--(c2);		 
		\draw[black, very thick, dashed] (c1) -- (c3);
		\draw [black, very thick, bend right=30] (c5) to (c1);
				}

	{}\hfill \hspace{-0.2cm} (a) $G^*$ \hspace{5.7 cm} (b) $H^*$\hfill {}			
 \caption{Graphs $G^*$ and $H^*$}
\label{f-G10}
\end{figure}

An example $(G^*,e,H^*,f)$ in $\Phi'$ was found
by Dr. Marion C. Gray~\cite{tutte} in 1930s, where $G^*$ and $H^*$ are graphs
shown in  Figure~\ref{f-G10}.
%mentioned by Tutte~\cite{tutte}
These two graphs are not isomorphic 
and even have non-isomorphic cyclic matroids,
because $H^*$, unlike $G^*$, contains 
a triangle having no common edge 
with any other triangle \cite{tutte}. 
However, $G^*\backslash e\cong H^*\backslash f$ and 
$G^*/e\cong H^*/f$,
where $e$ and $f$ are the edges in $G^*$ and $H^*$ 
which are displayed by dashed lines in Figure~\ref{f-G10}.
Thus $(G^*,e,H^*,f)\in \Phi'$.
Motivated by this example, we shall introduce a method in Section~\ref{sec5} to generate infinitely many members in $\Phi'$, as a partial answer to Problem~\ref{prob2}.

\section{Preliminary results
\label{sec1}}

In the following sections, we  restrict  to %the $T$-equivalent problem on
connected graphs due to the following fact indicated by (\ref{tuttep2}).

\begin{pro}\label{prop0}
If $G$ is disconnected with components $G_1,G_2,\dots, G_r$, 
or $G$ is connected with blocks $G_1,G_2,\dots, G_r$,
then 
\begin{equation}\label{prop0-eq1}
T_{G}(x,y)=\prod_{1\le i\le r}T_{G_i}(x,y).
\end{equation}
\end{pro}

Let $G=(V,E)$ be any graph. 
For any subset $V'$ of $V$, 
let $G\cdot V'$ denote the graph obtained from 
$G$ by identifying all vertices in $S$ as one vertex. 
Thus each edge in $G$ is still an edge in $G\cdot V'$,
although every edge of $G$ joining two vertices in $V'$ 
becomes a loop in $G\cdot V'$.
In particular, for any two distinct vertices $u,v$ in $G$,
we write $G\cdot uv$ for $G\cdot \{u,v\}$. 
%For any $e\in E$ with ends $u$ and $v$, let $G\backslash e$ be the graph obtained from $G$  by deleting $e$ and let $G/e$ denote the graph $(G\backslash e)\cdot uv$.We also say $G/e$ is obtained from $G$ by contracting $e$.
For any $A\subseteq E$,  
let $G\backslash A$ be the graph obtained from $G$ 
by deleting all edges in $A$.
Two edges $e_1$ and $e_2$ in $G$ are said to be \textit{parallel}
if they have the same pair of ends.
A relation among $T_{G}(x,y), T_{G\backslash e_1}(x,y)$
and $T_{G\backslash \{e_1,e_2\}}(x,y)$
is presented below.

\begin{pro}\label{prop1}
Assume that $e_1$ and $e_2$ are parallel edges in  $G$ 
which are not loops.
If $G\backslash \{e_1,e_2\}$ is connected, then 
\begin{equation}\label{prop1-eq1}
T_{G}(x,y)=(y+1)T_{G\backslash e_1}(x,y)
-y T_{G\backslash \{e_1,e_2\}}(x,y).
\end{equation}
\end{pro}

\proof 
By the given conditions, 
both $e_1$ and $e_2$ are not loops, 
$e_1$
is not a bridge of $G$
and $e_2$ is not a bridge of $G\backslash e_1$.

Assume that $u$ and $v$ are the two ends of $e_1$ and $e_2$, 
and let $F$ be the set of 
all parallel edges 
with ends $u$ and $v$, 
and let $k=|F|\ge 2$. 
Applying (\ref{tuttep2}) repeatedly yields that 
\begin{equation}\label{prop1-eq2}
	T_{G}(x,y)=T_{G\backslash e_1}(x,y)
	+y^{k-1}T_{(G\backslash F)\cdot uv}(x,y),
\end{equation}
and 
\begin{equation}\label{prop1-eq3}
	T_{G\backslash e_1}(x,y)=T_{G\backslash \{e_1,e_2\}}(x,y)
	+y^{k-2}T_{(G\backslash F)\cdot uv}(x,y).
\end{equation}
Thus, (\ref{prop1-eq1}) can be obtained 
from (\ref{prop1-eq2}) and (\ref{prop1-eq3}) by deleting 
$T_{(G\backslash F)\cdot uv}(x,y)$.
\proofend

%Let $H$ be any graph with parallel edges $f_1,f_2$ joining two distinct vertices. 
A consequence of Proposition~\ref{prop1} is then obtained. 

\begin{cor}\label{prop1-cor1}
Let $e_1,e_2$ be parallel edges in $G$ 
and $f_1,f_2$ be parallel edges in $H$. 
If both $G\backslash \{e_1,e_2\}$ and $H\backslash \{f_1,f_2\}$
are connected,
then any two of the 
equalities below imply another one:
\begin{equation}\label{prop1-cor2}
\left \{ 
\begin{array}{l}
T_{G\backslash \{e_1,e_2\}}(x,y)=T_{H\backslash \{f_1,f_2\}}(x,y), \\
T_{G\backslash e_1}(x,y)=T_{H\backslash f_1}(x,y),\\
T_{G}(x,y)=T_{H}(x,y).
\end{array}
\right.
\end{equation}
\end{cor}

\section{Necessary and sufficient conditions
\label{sec2}}

In this section, let $G$ and $H$ be two vertex-disjoint connected graphs.
Assume that  
$u_1,u_2,\dots,u_k$ are pairwise distinct vertices in $G$, and 
$v_1,v_2,\dots,v_k$
are pairwise distinct vertices in $H$.
We shall find a necessary and sufficient condition
for the following statement to be held:

\begin{quote}
(*) $G(u_1,\dots,u_k)\sqcup 
W(w_1,\dots,w_k)\tsim H(v_1,\dots,v_k)\sqcup 
W(w_1,\dots,w_k)$
for an arbitrary graph $W$,
where $w_1,w_2,\dots,w_k$ are
pairwise distinct vertices in $W$.
\end{quote}

For any $S\subseteq \{\{i,j\}: 1\le i<j\le k\}$,
let $G_S$ be the graph obtained from $G$ by adding 
a new edge joining $u_i$ and $u_j$ 
for each $\{i,j\}\in S$.
Observe that $G$ is a spanning subgraph of $G_S$
which has $|S|$ more edges than $G$\footnote{$G_S$ is a multiple graph.}.
Suppose that $H_S$ is defined similarly as $G_S$
for each $S\subseteq \{\{i,j\}: 1\le i<j\le k\}$.

Observe that for any $S\subseteq \{\{i,j\}: 1\le i<j\le k\}$,
$G_S$ is actually isomorphic to the graph 
$G(u_1,u_2,\dots,u_k)\sqcup 
W_0(w_1,w_2,\dots,w_k)$, where $W_0$ is the simple graph 
with vertex set $\{w_i: i\in [k]\}$ 
and edge set $\{w_iw_j: (i,j)\in S\}$.
Thus a necessary condition for statement (*)
is that $G_S\tsim H_S$
for every $S\subseteq \{\{i,j\}: 1\le i<j\le k\}$.
In the following, we shall show that this condition is also 
sufficient. 

\begin{theo}\label{them1}
If $G_S\tsim H_S$
%have the same Tutte polynomial 
for every subset $S$ of $\{\{i,j\}: 1\le i<j\le k\}$,
then for an arbitrary graph $W$ with 
pairwise
distinct vertices $w_1,w_2,\dots,w_k$, 
$G(u_1,u_2,\dots,u_k)\sqcup 
W(w_1,w_2,\dots,w_k)\tsim H(v_1,v_2,\dots,v_k)\sqcup 
W(w_1,w_2,\dots,w_k)$.
%have the same Tutte polynomial.
\end{theo}

\proof
Suppose that the result does not hold. 
Assume that $G,H,W$ are graphs stated above 
such that the result fails and $|V(W)|+|E(W)|$ has the 
minimum value. 
Let $Q=\{w_i: i\in [k]\}$.
By the assumption on $W$, 
$d_W(w)>0$ for each $w\in V(W)\setminus Q$.

We will complete the proof by showing the following claims.

\noindent {\bf Claim 1}: $W$ has no loops.

Suppose that $W$ contains a loop $e$. Then by (\ref{tuttep2}),
the Tutte polynomial of $G(u_1,u_2,\dots,u_k)\sqcup 
W(w_1,w_2,\dots,w_k)$ is equal to 
\begin{equation}\label{them1-eq4}
yT_{G(u_1,u_2,\dots,u_k)\sqcup 
W\backslash e(w_1,w_2,\dots,w_k)}(x,y).
\end{equation}
The Tutte polynomial of $H(v_1,v_2,\dots,v_k)\sqcup 
W(w_1,w_2,\dots,w_k)$ has a similar expression
as (\ref{them1-eq4}) where $G$ is changed to $H$
and each $u_i$ is changed to $v_i$.
By the minimality of $|V(W)|+|E(W)|$, 
the theorem holds for $W\backslash e$.
Thus it also holds for $W$, a contradiction.
Claim 1 holds.

\noindent {\bf Claim 2}: 
$V(W)=Q$.

Suppose the claim fails. 
Then $W$ has an edge $e$ with at most one end  in $Q$. 
Since $G$ and $H$ are connected, $e$ is a bridge in $G(u_1,u_2,\dots,u_k)\sqcup 
W(w_1,w_2,\dots,w_k)$ if and only if $e$ is a bridge in $H(v_1,v_2,\dots,v_k)\sqcup 
W(w_1,w_2,\dots,w_k)$.
Further in this case,
the Tutte polynomial of $G(u_1,u_2,\dots,u_k)\sqcup 
W(w_1,w_2,\dots,w_k)$ is equal to 
\begin{equation}\label{them1-eq1}
xT_{G(u_1,u_2,\dots,u_k)\sqcup 
W/e(w_1,w_2,\dots,w_k)}(x,y),
\end{equation}
 and otherwise, it is equal to 
\begin{equation}\label{them1-eq2}
T_{G(u_1,u_2,\dots,u_k)\sqcup 
W\backslash e(w_1,w_2,\dots,w_k)}(x,y)
+
T_{G(u_1,u_2,\dots,u_k)\sqcup 
W/e(w_1,w_2,\dots,w_k)}(x,y).
\end{equation}
Note that 
the Tutte polynomial of $H(v_1,v_2,\dots,v_k)\sqcup 
W(w_1,w_2,\dots,w_k)$
has a similar expression as 
(\ref{them1-eq1}) or (\ref{them1-eq2}), 
where $G$ is changed to $H$ and each $u_i$ is changed to $v_i$.

By assumption on the minimality of $|V(W)|+|E(W)|$, 
the theorem holds for both graphs $W/e$ and $W\backslash e$.
By the above results on the Tutte polynomials
of the two graphs $G(u_1,u_2,\dots,u_k)\sqcup 
W(w_1,w_2,\dots,w_k)$ 
and
$H(v_1,v_2,\dots,v_k)\sqcup 
W(w_1,w_2,\dots,w_k)$,
the theorem also holds for $W$, a contradiction. 

Hence Claim 2 holds.

\noindent {\bf Claim 3}: $W$ is a simple graph.

Suppose that $W$ has parallel edges $e_1$ and $e_2$.
Since $G$ and $H$ are connected and  $V(W)=Q$, $G(u_1,u_2,\dots,u_k)\sqcup 
W\backslash \{e_1,e_2\}(w_1,w_2,\dots,w_k)$ and  $H(v_1,v_2,\dots,v_k)\sqcup 
W\backslash \{e_1,e_2\}(w_1,w_2,\dots,w_k)$ are connected.
By the assumption on the minimality of $|V(W)|+|E(W)|$,
the theorem holds for both $W\backslash e_1$ and 
$W\backslash \{e_1,e_2\}$.
By Corollary~\ref{prop1-cor1}, 
the theorem also holds for $W$, 
a contradiction. 
Claim 3 holds.

\noindent {\bf Claim 4}: The result holds for $G,H$ and $W$.
 
By Claims 1--3,  $G(u_1,u_2,\dots,u_k)\sqcup 
W(w_1,w_2,\dots,w_k)$ 
and $H(v_1,v_2,\dots,v_k)\sqcup 
W(w_1,w_2,\dots,w_k)$ 
are actually respectively 
the graphs $G_S$ and $H_S$ 
for some $S\subseteq \{\{i,j\}: 1\le i<j\le k\}$.
Thus Claim 4 holds by the given condition.

Claim 4 contradicts the given assumption. The result is proven.
\proofend

\iffalse
For any partition $P$ of $\{1,2,\dots,k\}$ 
into non-empty subsets $I_1,\dots,I_r$, 
let $G(P)$ denote the graph obtained from $G$ by identifying all vertices 
in $\{u_j: j\in I_i\}$ to give a single vertex for all $i=1,2,\dots, r$.
$H(P)$ is defined similarly.
%Let ${\cal P}$ denote the set of all such partitions $P$.

\begin{theo}\label{them1-1}
	If $G(P)$ and $H(P)$ are $T$-equivalent
	%have the same Tutte polynomial 
	for every partition $P$ 
	of $\{1,2,\dots,k\}$, then 
	for an arbitrary graph $W$ with 
	distinct vertices $w_1,w_2,\dots,w_k$, 
	$G(u_1,u_2,\dots,u_k)\sqcup 
	W(w_1,w_2,\dots,w_k)$
	and $H(v_1,v_2,\dots,v_k)\sqcup 
	W(w_1,w_2,\dots,w_k)$ are $T$-equivalent. 
\end{theo}
\fi

A necessary condition for Theorem~\ref{them1}
is quite obvious. 

\begin{pro}\label{nece-con}
The conditions in Theorem~\ref{them1} imply that 
$G$ and $H$ have the same number of loops 
and,  for each pair of numbers $i,j: 1\le i<j\le k$, 
the number of edges in $G$ joining $u_i$ and $u_j$ 
is equal to the 
number of edges in $H$ joining $v_i$ and $v_j$.
\end{pro}

\proof 
Our proof will apply the basic fact that 
for any graph $G$, the maximum integer $r$ such that $y^r$ is a factor of $T_G(x,y)$ is the number of loops in $G$. Then $G$ and $H$ must have the same number of loops as $T_{G_S}(x,y)=T_{H_S}(x,y)$ when $S=\emptyset$ by the given condition.

Also, the condition that $T_{G_S}(x,y)=T_{H_S}(x,y)$ when $S=\{\{i,j\}\}$ indicates that $T_{G_S/ u_iu_j}(x,y)=T_{H_S/v_iv_j}(x,y)$ for each pair of numbers $i,j: 1\le i<j\le k$. Therefore, graphs $G_S/u_iu_j$ and $H_S/v_iv_j$ have the same number of loops, which further implies that the number of edges in $G$ joining $u_i$ and $u_j$ 
is equal to the number of edges in $H$ joining $v_i$ and $v_j$.
\proofend

We end this section by providing another necessary 
and sufficient condition for statement (*) to be held.
%condition for Theorem~\ref{them1}.

%\section{Alternate conditions for Theorem~\ref{them1}}

\def \setp{{\cal P}}
\def \bkp{{\bf P}}

Let $\setp_k$ be the set of partitions of $[k]$.
For any partition $\bkp=\{P_i: i\in [r]\}\in \setp_k$,
let $G(\bkp)$ be the graph 
%$\cdots((G\cdot U_1)\cdot U_2)\cdots) \cdot U_r$,where  $U_j=\{u_i: i\in P_j\}$ for $j=1,2,\dots,r$.
%Thus $G(P)$ is the graph 
obtained from $G$ 
by identifying all vertices in the $\{u_i: i\in P_j\}$
as one vertex $x_j$ for all $j\in [r]$.
Assume that 
$H(\bkp)$ is the graph defined similarly as $G(\bkp)$.
%i.e., $H(P)$ is the graph obtained from $H$ by identifying all vertices in $\{v_i: i\in P_j\}$ as one vertex $y_j$ for all $j=1,2,\dots,r$.
It is not difficult to prove 
by applying Theorem~\ref{them1}
and (\ref{tuttep2}) that 
 statement (*) implies that 
$G(\bkp)\tsim H(\bkp)$
for every $\bkp\in \setp_k$.
In fact, this condition is sufficient as well.

\begin{theo}\label{them2}
%Assume that $\{u_1,u_2,\dots,u_k\}$ is an independent set of $G$ and $\{v_1,v_2,\dots,v_r\}$ is an independent set of $H$.
%Assume that both $G$ and $H$ are connected graphs with $\{u_1,u_2,\dots,u_k\}\subseteq V(G)$ and $\{v_1,v_2,\dots,v_k\}\subseteq V(H)$.
If $G(\bkp)\tsim H(\bkp)$
%have the same Tutte polynomial 
for every $\bkp\in \setp_k$, then 
$G(u_1,\dots,u_k)\sqcup 
W(w_1,\dots,w_k)\tsim H(v_1,\dots,v_k)\sqcup 
W(w_1,\dots,w_k)$
for an arbitrary graph $W$, where 
$w_1,w_2,\dots,w_k$ 
are pairwise distinct vertices in $W$.
\end{theo}

\proof Assume that $G(\bkp)\tsim H(\bkp)$ for each $\bkp\in \setp_k$.
By Theorem~\ref{them1},
it suffices to show that 
for each subset $S$ of $\{\{s,t\}: 1\le s<t\le k\}$, $G_S\tsim H_S$ holds.
%have the same Tutte polynomial.

Let $S$ be any subset of $\{\{s,t\}: 1\le s<t\le k\}$.
Write $E(G_S)-E(G)=\{e_{s,t}: \{s,t\}\in S\}$,
where $e_{s,t}$ is an edge joining $u_s$ and $u_t$. 
As $G$ is connected, 
applying the last two equalities in (\ref{tuttep2}) repeatedly 
%(i.e., the one on deleting and contracting an edgewhich is neither a loop nor a bridge) 
on all edges in the set $\{e_{s,t}: \{s,t\}\in S\}$,
%where $e_{i,j}$ is an edge joining $u_i$ and $u_j$,
%until there is no such edge, 
$T_{G_S}(x,y)$ has the following expression:
\begin{equation}\label{them2-eq1}
T_{G_S}(x,y)=\sum_{\bkp\in \setp_k} a(\bkp,S)y^{n(\bkp,S)}T_{G(\bkp)}(x,y),
\end{equation}
where
$a(\bkp,S)$ and $n(\bkp,S)$ are non-negative integers 
depending on $\bkp$ and $S$ only. 
Actually, for any $\bkp=\{P_i: i\in [r]\}\in \setp_k$,
$a(\bkp,S)$ is the number of those spanning forests $F$
of the graph $N_S$, 
where $N_S$ is the simple graph with vertex set 
$\{u_i: i\in [k]\}$
and edge set $\{u_su_t: \{s,t\}\in S\}$,
such that 
$F$ has exactly $r$ components whose vertex sets 
are $\{u_i: i\in P_j\}$ for $j\in [r]$,
and $n(\bkp,S)=e(\bkp,S)-k+r$, 
where $e(\bkp,S)$ is the size of 
the set 
$\{\{s,t\}\in S: s,t\in P_j, j\in [r]\}$.

Thus we also have 
\begin{equation}\label{them2-eq2}
T_{H_S}(x,y)=\sum_{\bkp\in \setp_k}a(\bkp,S) y^{n(\bkp,S)}T_{H(\bkp)}(x,y).
\end{equation}
Hence $T_{G_S}(x,y)=T_{H_S}(x,y)$ holds for 
every subset $S$ of $\{\{i,j\}: 1\le i<j\le k\}$.
By Theorem~\ref{them1}, the result holds.
\proofend

Theorem~\ref{them2} is important for proving the 
main result in the next section.

\section{Extension of 
	%Tutte's work in 
	Theorem~\ref{theo-tutte}
\label{secmain}}

In this section, assume that 
$R$ is a connected  
graph with 
$\{u_i: i\in [k]\}\subseteq V(R)$ 
and an automorphism $\psi$
such that $\psi(u_i)=u_{i+1}$ for all $i\in [k]$,
where $u_{i+k}$ and $u_i$ always denote the same vertex in  $R$. 
We shall extend %Tutte's result in 
Theorem~\ref{theo-tutte}
and 
show that for some integers $k\ge 6$,  $R(u_1,\dots,u_{k})\sqcup W(w_1,\dots,w_{k})$
and $R(u_{k},\dots,u_1)\sqcup W(w_1,\dots,w_{k})$ are $T$-equivalent
for any graph $W$ satisfying certain conditions, where 
$w_1,\dots,w_{k}$ are distinct vertices in $W$.

First we establish the following 
preliminary result.
%for proving our main result in this section. 

\begin{pro}\label{isomor}
Let $W$ be a loopless graph and
$w_1,w_2,\dots,w_k$ be $k$ pairwise 
distinct vertices in $W$.
If $W$ has an automorphism $\rho$ 
such that 
$\rho(w_{1+s})=w_{a-s}$
for all integers $s$ with $0\le s\le k-1$,
where $a\in [k]$
%$a$ is an integer with $1\le a\le k$ 
and $w_{s_1}$ and $w_{s_2}$ 
represent the same vertex in $W$
whenever $s_1-s_2\equiv 0 \pmod{k}$, 
then $R(u_1,\dots,u_{k})\sqcup W(w_1,\dots,w_{k})$ 
and $R(u_{k},\dots,u_1)\sqcup W(w_1,\dots,w_k)$
are isomorphic.
\end{pro}

\proof 
Let $G_1,G_2,G_3$ and $G_4$ denote 
the graphs in statements (i) to (vi) below respectively:
\begin{enumerate}
\item $R(u_{1},\dots,u_k)\sqcup W(w_k,\dots,w_1)$
(i.e., $R(u_{k},\dots,u_1)\sqcup W(w_1,\dots,w_k)$);
\item $R(u_{1-a},u_{2-a},\dots,u_{k-a})
\sqcup W(w_k,w_{k-1},\dots,w_1)$;
%\item $R(u_{a+1},u_{a+2},\dots,u_{a+k})\sqcup W(w_k,w_{k-1},\dots,w_1)$;

\item $R(u_{1},u_2,\dots,u_{k})
\sqcup W(w_{a},w_{a-1},\dots,w_1, w_k,\dots, w_{a+1})$;

\item $R(u_1,u_2,\dots,u_{k})\sqcup W(w_1,w_2,\dots,w_{k})$.
%where  $W$ has an automorphism $\rho$ such that $\rho(w_{1+s})=w_{a-s}$ for all $s=0,1,\dots,k-1$.
\end{enumerate}
It can be verified that 
$G_1\cong G_2\cong G_3\cong G_4$.
It is obvious that $G_2\cong G_3$. 
The isomorphism of $G_1$ and $G_2$ 
follows from the fact that 
$\psi^{k-a}(u_i)=u_{i-a}$ for all $i\in [k]$,
as $\psi$ is an automorphism of $R$ 
such that $\psi(u_i)=u_{i+1}$ 
for all $i\in [k]$,
while the isomorphism of $G_3$ and $G_4$ 
follows from the condition 
that $W$ has an automorphism $\rho$
such that $\rho(w_{1+s})=w_{a-s}$ holds 
for all integers $s$ with $0\le s\le k-1$.

Thus the result holds.
\proofend

\noindent\textbf{Remark:} Proposition~\ref{isomor} also holds 
if the automorphism $\psi$ of $R$ is not provided 
but the condition for $\rho$ holds for $a=k$.

%Now assume that $k=rg$, where $r$ and $g$ are positive integers. 
Let $W_0$ be a loopless graph 
with disjoint subsets 
$V_1$ and $V_2$ of vertices, where 
$V_1=\{w_i: i\in [rg]\}$
%w_1,w_2,\dots, w_{rg}\}$
and $V_2=\{x_i: i\in [r]\}$,
where $w_{j+rg}$ and $w_j$ (resp. $x_{j}$ and $x_{j+r}$)
denote the same vertex in $W_0$
for all integers $j$. 
It is possible that $V_1\cup V_2$
is a proper subset of $V(W_0)$.
Assume that $W_0$ has two automorphisms $\phi$ and $\rho$ 
satisfying the following conditions:
\begin{enumerate}
\item  
$\phi(w_i)=w_{i+g}$ for all $i\in [rg]$
%for all $i=1,2,\dots,rg$,
and $\phi(x_s)=x_{s+1}$ for all 
$i\in [r]$;
%for all $s=1,\dots,r$;
\item for some $c\in [rg]$, 
%some integer $c$ with $1\le c\le rg$, 
$\rho(w_{1+s})=w_{c-s}$ holds
for all integers $s$ with $0\le s\le rg-1$, and for some $c'\in [r]$, 
$\rho(x_{1+i})=x_{c'-i}$ holds
for all integers $i$ with $0\le i\le r-1$.
\end{enumerate}
Two examples for such graphs $W_0$ 
are shown in Figure~\ref{f-G11}.

\tikzstyle{cblue}=[circle, draw, thin,fill=blue!20, scale=0.5]
\begin{figure}[!ht]
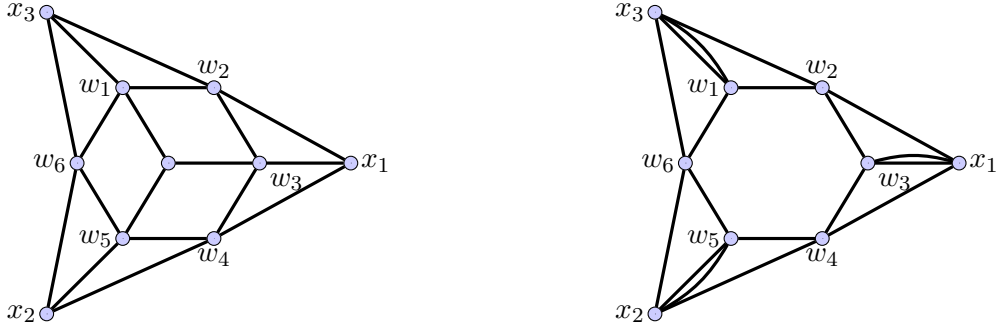

	\tikzp{1}
	{
		\foreach \place/\y in {{(0,0)/1}, {(0.6,1)/2},{(1.2,0)/3}, {(0.6,-1)/4}, {(-0.6,-1)/5}, {(-1.2,0)/6},{(-0.6,1)/7},{(-1.6,2)/8}, {(2.4,0)/9},{(-1.6,-2)/10}}   
		\node[cblue] (b\y) at \place {};
		
		\filldraw[black] (b1) circle (0pt)node[anchor=east] {};
		\filldraw[black] (b2) circle (0pt)node[anchor=south] {$w_2$};
		\filldraw[black] (b3) circle (0pt)node[anchor=north west] {$w_3$};
		\filldraw[black] (b4) circle (0pt)node[anchor=north] {$w_4$};
		\filldraw[black] (b5) circle (0pt)node[anchor=east] {$w_5$};
		\filldraw[black] (b6) circle (0pt)node[anchor=east] {$w_6$};
		\filldraw[black] (b7) circle (0pt)node[anchor=east] {$w_1$};
		\filldraw[black] (b8) circle (0pt)node[anchor=east] {$x_3$};
		\filldraw[black] (b9) circle (0pt)node[anchor=west] {$x_1$};
		\filldraw[black] (b10) circle (0pt)node[anchor=east] {$x_2$};
			
		\draw[black, very thick] (b1) -- (b3) -- (b4) -- (b5) --(b1) --(b7)--(b6)--(b5) -- (b10);		 
		\draw[black, very thick] (b7) -- (b2) -- (b3) -- (b9);
		\draw[black, very thick] (b7) -- (b8) -- (b2) --(b9) -- (b4) -- (b10) -- (b6) -- (b8);
				
		\foreach \place/\y in { {(0.6+8,1)/2},{(1.2+8,0)/3}, {(0.6+8,-1)/4}, {(-0.6+8,-1)/5}, {(-1.2+8,0)/6},{(-0.6+8,1)/7},{(-1.6+8,2)/8}, {(2.4+8,0)/9},{(-1.6+8,-2)/10}}   
		\node[cblue] (c\y) at \place {};

		\filldraw[black] (c2) circle (0pt)node[anchor=south] {$w_2$};
		\filldraw[black] (c3) circle (0pt)node[anchor=north west] {$w_3$};
		\filldraw[black] (c4) circle (0pt)node[anchor=north] {$w_4$};
		\filldraw[black] (c5) circle (0pt)node[anchor=east] {$w_5$};
		\filldraw[black] (c6) circle (0pt)node[anchor=east] {$w_6$};
		\filldraw[black] (c7) circle (0pt)node[anchor=east] {$w_1$};
		\filldraw[black] (c8) circle (0pt)node[anchor=east] {$x_3$};
		\filldraw[black] (c9) circle (0pt)node[anchor=west] {$x_1$};
		\filldraw[black] (c10) circle (0pt)node[anchor=east] {$x_2$};
			
		\draw[black, very thick] (c3) -- (c4) -- (c5);
		\draw[black, very thick] (c7)--(c6)--(c5) -- (c10);		 
		\draw[black, very thick] (c7) -- (c2) -- (c3) -- (c9);
		\draw[black, very thick] (c7) -- (c8) -- (c2) --(c9) -- (c4) -- (c10) -- (c6) -- (c8);

		\draw [black, very thick, bend right=15] (c7) to (c8);
		\draw [black, very thick, bend left=15] (c3) to (c9);
		\draw [black, very thick, bend left=15] (c5) to (c10);
		
}
		
 \caption{Examples for graphs $W_0$ with $r=3$ and $g=2$} 
\label{f-G11}
\end{figure}

We are now ready 
to prove the main result in this section
by applying Theorem~\ref{them2} and 
Proposition~\ref{isomor}.

%\newpage 

\begin{theo}\label{more-rotors}
%Let $r$ and $g$ be integers with $1\le r\le 5$ and $g\ge 1$.

Let $Y$ be an arbitrary graph with a subset of vertices 
$y_1,\dots,y_{r}$, and let 
$W$ denote the graph $W_0(x_1,\dots,x_{r})\sqcup 
Y(y_1,\dots,y_{r})$.
\iffalse
Let $W$ denote the graph $W_0(x_1,\dots,x_{r}, z_1,\dots,z_h)\sqcup 
Y(y_1,\dots,y_{r+h})$, where
$z_1,\dots,z_h$ are fixed vertices of both $\phi$ and $\rho$ in $W_0$ and $Y$ 
is 
an arbitrary graph with a subset of vertices 
$y_1,\dots,y_{r+h}$.
\fi
If $r\le 5$
%, $g\ge 1$
and $k=rg$,
then $R(u_1,\dots,u_{rg})\sqcup W(w_1,\dots,w_{rg})$ and $ R(u_{rg},\dots,u_1)\sqcup W(w_1,\dots,w_{rg})$
are $T$-equivalent.
\end{theo}

\def \setq {{\cal Q}}

\proof 
For the sake of convenience, 
let $L$ and $F$ 
denote graphs $R(u_1,\dots,u_{rg})\sqcup W_0(w_1,\dots,w_{rg})$ 
and $R(u_{rg},\dots,u_1)\sqcup W_0(w_1,\dots,w_{rg})$
respectively, and 
let $\Gamma$ be the family of 
partitions of the set 
$B=\{x_1,\dots,x_{r}\}$.
%$B=\{x_1,\dots,x_{r}, z_1,\dots,z_h\}$.
For any $\setq=
\{Q_i: i\in [t]\}\in \Gamma$, 
%\{Q_1,Q_2,\dots,Q_t\}\in \Gamma$, 
let $L(\setq)$ denote the graph obtained from $L$ 
by identifying all vertices in $Q_j$ as one vertex for each 
$j\in [t]$, and $F(\setq), W_0(\setq)$ are defined similarly.
By Theorem~\ref{them2},
it suffices to show that $L(\setq)$ and $F(\setq)$ have the 
same Tutte polynomial for every $\setq\in \Gamma$.

Now let S1, S2, S3 and S4 denote the following statements:
\begin{enumerate}
\item [S1:] for every $\setq\in \Gamma$, 
$L(\setq)$ and $F(\setq)$ have the 
same Tutte polynomial;
\item[S2:] for every $\setq\in \Gamma$,
$L(\setq)\cong F(\setq)$ holds;
\item[S3:] for every $\setq\in \Gamma$, 
$W_0(\setq)$ has an automorphism $\eta$ 
and an integer $a'$ with $1\le a'\le rg$
such that 
$\eta(w_{1+s})=w_{a'-s}$ for all $s=0,1,\dots,rg-1$.
\item[S4:]
for every $\setq=\{Q_i: i\in [t]\}\in \Gamma$, 
there exist integers $p$ and $d$, where $p\in [rg]$,
such that the two conditions below are satisfied
with $\pi=\rho \phi^d$:
\begin{enumerate}
\item [(a)] $\pi(w_{1+s})=w_{p-s}$ 
for all integers $s$ with $0\le s\le rg-1$; and 
\item [(b)] $\pi(\setq)=
\{\pi(Q_i): i\in [t]\}=\setq$.
\end{enumerate}
\end{enumerate}
It can be proved that 
the relation 
``S1 $\Leftarrow$ S2 $\Leftarrow$ S3 $\Leftarrow$ S4" holds,
where ``S1 $\Leftarrow$ S2" means that 
$S1$ follows from $S2$.
The relation ``S1 $\Leftarrow$ S2" is trivial 
and ``S2 $\Leftarrow$ S3" holds by 
Proposition~\ref{isomor}. 
The relation ``S3 $\Leftarrow$ S4" holds 
as an automorphism $\eta$ of $W_0(\setq)$ in S3 can be 
obtained from the automorphism $\pi$ of $W_0$ in S4
by the permutation $\delta$ of $1,2,\dots,t$ 
such that $\pi(Q_i)=Q_{\delta(i)}$ holds 
for all $i\in [t]$.

Thus it remains to show that S4 holds.
%As each vertex $z_j$ is a fixed vertex of $\rho$ and $\phi$, $z_j$ is also a fixed vertex of $\rho\phi^d$ for any integer $d$. Hence it suffices to show that  S4 holds for the case that $B=\{x_1,\dots,x_r\}$.
%We now assume that $\Gamma$ is the family of partitions of $B=\{x_1,\dots,x_r\}$.
For any integer $b,b'$ with $1\le b,b'\le r$, 
let $\pi_{(b,b')}$ denote the automorphism 
$\rho\phi^d$ of $W_0$, where $d=c'+1-b-b'$. 
For any integer $s$ with $0\le s\le rg-1$, 
\begin{equation}\label{more-rotors-eq1}
\pi_{(b,b')}(w_{1+s})=\rho( \phi^d(w_{1+s}))
=\rho(w_{1+s+dg})
=w_{c-dg-s}=w_{p-s},
\end{equation}
where $p$ is the integer with $1\le p\le rg$ 
and $p\equiv c-dg\ (\mbox{mod\ }rg)$, and 
for any integer $i$, 
\begin{equation}\label{more-rotors-eq2}
\pi_{(b,b')}(x_{b+i})=\rho( \phi^d(x_{b+i}))
=\rho(x_{b+i+d})
=x_{c'-(b+i+d-1)}=x_{b'-i}.
\end{equation}
By (\ref{more-rotors-eq1}),
condition (a) in S4 is satisfied 
with $\pi_{(b,b')}$ for any integers $b$ and $b'$ 
with $1\le b,b'\le r$.
In the following, 
it suffices to show that 
for every $\setq\in \Gamma$, there exist  
integers $b$ and $b'$ 
with $1\le b,b'\le r$ such that $\pi_{(b,b')}(\setq)=\setq$.

Write $\pi_b$ for $\pi_{(b,b)}$.
We complete the proof by 
considering all possible 
$\setq=\{Q_i: i\in [t]\}\in \Gamma$, where 
$|Q_1|\ge \cdots\ge |Q_t|$.
\begin{enumerate}
\item[Case 1:] $|Q_1|=1$ or $|Q_1|=r$. 
%Let $\pi=\rho$. 
It is obvious that %$\rho=\pi_{(1,r)}$ and 
$\pi_{(1,r)}(\setq)=\setq$.

\item[Case 2:]  $|Q_1|=r-1$.
Then there exists an integer $b$ with $1\le b\le r$ 
such that $Q_2=\{x_b\}$.
% and $Q_1=B-Q_2$.
By (\ref{more-rotors-eq2}), we have 
$\pi_b(Q_i)=Q_i$ for $i=1,2$ 
and so $\pi_b(\setq)=\setq$.

\item[Case 3:]  $|Q_1|=r-2$.
Then there exist integers $b,b'$ with $1\le b,b'\le r$ 
such that either $Q_{2}=\{x_b\}$ and $Q_3=\{x_{b'}\}$, or $Q_2=\{x_b,x_{b'}\}$.
By (\ref{more-rotors-eq2}), we have 
$\pi_{b,b'}(\setq)=\setq$.

Note that by Cases 1--3, 
we now can assume that 
$r=5$ and $|Q_1|=2$.

\iffalse
\item[Case 3:] $r=4$ and $|Q_1|=2$.
Then, there are integers $b$ with $1\le b\le 4$ 
and $j$ with $j\in \{1,2\}$
such that $Q_1=\{x_b,x_{b+j}\}$
and either $Q_2=\{x_{b-1},x_{b-j-1}\}$ 
or $\{Q_2, Q_3\}=\{\{x_{b-1}\},\{x_{b-j-1}\}\}$. 
By (\ref{more-rotors-eq2}), 
%$\pi_{(b,b+j)}(b+i)=b+j-i$ holds for all integer $i$.
$\pi_{(b,b+j)}(Q_1)=Q_1$, 
$\pi_{(b,b+j)}(x_{b-1})=x_{b+j+1}$ and 
$\pi_{(b,b+j)}(x_{b-j-1})=x_{b+2j+1}$.
It is easy to verify that 
$\pi_{(b,b+j)}(\setq)=\setq$ holds. % for $j=1,2$.

\item[Case 4:]  $r=5$ and $|Q_1|=3$.
There are integers $b$ and $j$ 
with $1\le b\le 5$ and $j\in \{1,2\}$ 
such that $Q_1=\{x_{b-j},x_b,x_{b+j}\}$
and either $Q_2=\{x_{b-3+j},x_{b+3-j}\}$
or $\{Q_2,Q_3\}=\{\{x_{b-3+j}\},\{x_{b+3-j}\}\}$.
By (\ref{more-rotors-eq2}), we have 
$\pi_b(x_{b+i})=x_{b-i}$ for all $i\in \{0,-j,j,3-j,j-3\}$.
Thus  $\pi_b(\setq)=\setq$.
\fi

\item[Case 4:] $r=5$, $|Q_1|=2$ and $t=4$.
There are integers $b,b'$ with $1\le b,b'\le 5$,
such that $Q_1=\{x_{b},x_{b'}\}$ and 
$|Q_2|=|Q_3|=|Q_4|=1$.
By (\ref{more-rotors-eq2}), we have $\pi_{b,b'}(\setq)=\setq$.

\item[Case 5:] $r=5$, $|Q_1|=2$ and $t=3$.
There is an integer $b$ with $1\le b\le 5$ such that 
$Q_3=\{x_b\}$ and either 
$\{Q_1,Q_2\}
=\{\{x_{b-1},x_{b+1}\},\{x_{b-2},x_{b+2}\}\}$
or  
$\{Q_1,Q_2\}=\{\{x_{b-1},x_{b-1-j}\}, 
\{x_{b+1},x_{b+1+j}\}\}$ 
for some $j\in \{1,2\}$.
By (\ref{more-rotors-eq2}), we have 
$\pi_b(x_{b+i})=x_{b-i}$ for all $i\in \{1,-1,1+j,-1-j, 2,-2\}$.
Thus  $\pi_b(\setq)=\setq$.
\end{enumerate}
Thus we complete the proof.
\proofend

\noindent \textbf{Remark:} for any positive integers 
$r$ and $g$, 
if $W=W_0(x_1,\dots,x_{r})\sqcup 
Y(y_1,\dots,y_r)$
and $W'=W_0(x_r,\dots,x_{1})\sqcup 
Y(y_1,\dots,y_r)$, where $Y$ is 
an arbitrary graph with a subset of vertices 
$y_1,y_2,\dots,y_r$,
then $R(u_1,\dots,u_{rg})\sqcup W(w_1,\dots,w_{rg})$ 
and $R(u_{rg},\dots,u_{1})\sqcup W'(w_{1},\dots,w_{rg})$
are isomorphic and hence 
are $T$-equivalent.

%\section{Graphs satisfying conditions in Theorem~\ref{them1}}

Theorem~\ref{them1} provides a method of generating new 
$T$-equivalent graphs from two graphs 
$G$ and $H$ with $\{u_i: i\in [k]\}\subseteq V(G)$ 
and $\{v_i: i\in [k]\}\subseteq V(H)$
which satisfy the condition that $G_S\tsim H_S$ 
for every $S\subseteq \{\{i,j\}:1\le i<j\le k\}$.
What graphs $G$ and $H$ satisfy these conditions?
We are motivated to propose the following problem.

\begin{prob}\label{prob1}
Searching for graphs $G$ and $H$ with
$\{u_i: i\in [k]\}\subseteq V(G)$ and 
$\{v_i: i\in [k]\}\subseteq V(H)$
such that $G_S\tsim H_S$ 
for every $S\subseteq \{\{i,j\}:1\le i<j\le k\}$.
\end{prob}

If $k\in \{0,1\}$, the problem is just on 
searching for a pair of $T$-equivalent graphs $G$ and $H$.
It is trivial if $G$ is isomorphic to $H$ 
with $\{u_i: i\in [k]\}$ corresponding to
$\{v_i: i\in [k]\}$.

\section{Generating $T$-equivalent pairs
\label{sec5}}

\iffalse
The graph considered here are multigraphs. 
For any graph $G$, the Tutte polynomial of $G$ 
is defined as follows:
\begin{equation}\label{tuttep1}
T_G(x,y)=\sum_{A\subseteq E}(x-1)^{r(E)-r(A)}(y-1)^{|A|-r(A)},
\end{equation}
For any two graphs $G$ and $H$, 
we say $G$ and $H$ are {\it $T$-equivalent}, 
written as  $G\tsim H$, if $T_G(x,y)=T_H(x,y)$ holds. 

For some special cases of Tutte polynomials, 
such as reliability polynomials
and chromatic polynomials of graphs, 
Theorem~\ref{them1} can also be applied to them.
In other words, the result of Theorem~\ref{them1} also holds 
if ``two graphs are $T$-equivalent" is 
replaced by ``two graphs have the same chromatic polynomial"
or ``two graphs have the same reliability polynomial".
Thus Problems~\ref{prob1} and~\ref{prob2}
can also be changed to ones for chromatic polynomials
or reliability polynomials or any special case of Tutte polynomials. 

In this section, we present a method 
for generating members in $\Phi$
from a known member in $\Phi$.

we present some pairs of graphs 
$G$ and $H$ which are equivalent for some polynomial
such that $G\backslash e$ and $H\backslash f$ are also 
equivalent for this polynomial,
where $e\in E(G)$ and $f\in E(H)$.
It seems quite difficult to 
find such pairs of graphs $G$ and $H$ with their edges
$e$ and $f$.%  satisfying such conditions.

It will be seen later 
that $(G_0,e,H_0,f)$ is one of the members in 
$\Phi'$ generated by this 
method, % from $(G_0,e,H_0,f)\in \Phi'$,
where $G'_0$ and $H'_0$ are the graphs 
shown in Figure~\ref{f-G10-2-0}.
\fi

In this section, we present a method 
for generating members in $\Phi'$
from any known member $(G,e,H,f)$ in $\Phi'$ according to an isomorphism between $G/e$ and $H/f$. The specific steps are as follows.

\begin{enumerate}
\item[Step 1.]
%Since $G\backslash e\cong H\backslash f$, 
Vertices in $G$ and $H$ 
are labeled as $u_1,u_2,\dots,u_n$
and  $v_1,v_2,\dots,v_n$ respectively 
such that the mapping $\phi$ defined by $\phi(u_i)=v_i$ 
for all $i\in [n]$ is an isomorphism 
from $G\backslash e$ to $H\backslash f$, where $e$ and $f$ are edges in 
$G$ and $H$ respectively
such that $(G,e,H,f)$ in $\Phi'$ .
An example %for $(G_0,e,H_0,f)\in \Phi'$
 is shown in Figure~\ref{f-G10-2-0}.
%Note that more examples can be obtained by edge gluing a fixed edge of a graph copy to each edge in $G_0$ and $H_0$.

\tikzstyle{cblue}=[circle, draw, thin,fill=blue!20, scale=0.5]
\begin{figure}[!ht]
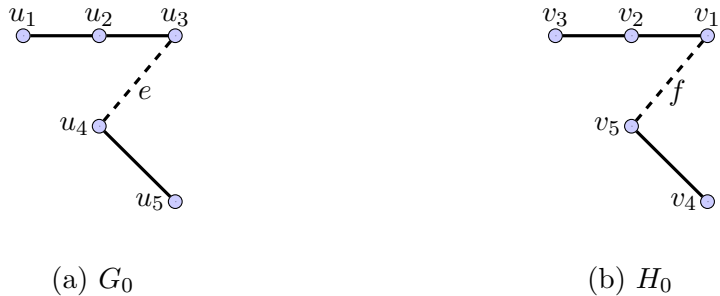

	\tikzp{1}
	{
		\foreach \place/\y in {{(0,0)/1}, {(0,1.2)/2},{(-1,1.2)/3},{(1,1.2)/4},{(1,-1)/5}}   
		\node[cblue] (b\y) at \place {};
		
		\filldraw[black] (b1) circle (0pt)node[anchor=east] {$u_4$};
		\filldraw[black] (b2) circle (0pt)node[anchor=south] {$u_2$};
		\filldraw[black] (b3) circle (0pt)node[anchor=south] {$u_1$};
		\filldraw[black] (b4) circle (0pt)node[anchor=south] {$u_3$};
		\filldraw[black] (b5) circle (0pt)node[anchor=east] {$u_5$};		
				\node [style=none] (cap1) at (0.6,0.45)  {$e$};
				
		\draw[black, very thick] (b3) -- (b2)--(b4);
		\draw[black, dashed, very thick] (b4) --(b1);
		\draw[black, very thick] (b5) --(b1);
		
		\foreach \place/\y in {{(0+7,0)/1}, {(0+7,1.2)/2},{(-1+7,1.2)/3},{(1+7,1.2)/4},{(1+7,-1)/5}}   
		\node[cblue] (c\y) at \place {};
		
		\filldraw[black] (c1) circle (0pt)node[anchor=east] {$v_5$};
		\filldraw[black] (c2) circle (0pt)node[anchor=south] {$v_2$};
		\filldraw[black] (c3) circle (0pt)node[anchor=south] {$v_3$};
		\filldraw[black] (c4) circle (0pt)node[anchor=south] {$v_1$};
		\filldraw[black] (c5) circle (0pt)node[anchor=east] {$v_4$};		
				\node [style=none] (cap1) at (0.6+7,0.45)  {$f$};
				
		\draw[black, very thick] (c3) -- (c2)--(c4);
		\draw[black, dashed, very thick] (c4) --(c1);
		\draw[black, very thick] (c5) --(c1);
									}

	{}\hfill \hspace{-0.2cm} (a) $G_0$ \hspace{5.7 cm} (b) $H_0$\hfill {}			
 \caption{$(G_0, e, H_0, f)\in \Phi'$}
\label{f-G10-2-0}
\end{figure}

\item[Step 2.]
Let $u_{s_1}$ and $u_{s_2}$ be the two ends of $e$ in $G$,
and $v_{t_1}$ and $v_{t_2}$ be the two ends of $f$ in $G$.
Let $u\idf{s_1}{s_2}$
%$\{u_{s_1},u_{s_2}\}$ 
denote the new vertex 
in $G/e$ obtained by identifying 
$u_{s_1}$ and $u_{s_2}$ in $G\backslash e$.
Thus, the vertex set of $G/e$ is
$$
%(\{u_1,u_2,\dots,u_n\}\setminus \{u_{s_1},u_{s_2}\})
\{u_i: i\in [n]\setminus \{s_1,s_2\}\}
\cup \{u\idf{s_1}{s_2}\}. 
$$ 
Similarly, the vertex set of $H/f$ is
$$
\{v_i: i\in [n]\setminus \{t_1,t_2\}\}
\cup \{v\idf{t_1}{t_2}\}.
%(\{v_1,v_2,\dots,v_n\}\setminus \{v_{t_1},v_{t_2}\})\cup \{v\idf{t_1}{t_2}\},
$$ 
%where $v\idf{t_1}{t_2}$ denote the new vertex in $H/f$ obtained by identifying $v_{t_1}$ and $v_{t_2}$ in $H\backslash f$.

Let $\psi$ be an isomorphism from $G/e$ to $H/f$.
For example, for the two graphs in Figure~\ref{f-G10-2-0},
we may either choose $\psi$ such that  
\begin{equation}\label{sec5-eq1}
(\psi(u_1),\psi(u_2),
\psi(u\idf{3}{4}),\psi(u_5)) 
=(v_4, v\idf{1}{5},v_2,v_3)
\end{equation}
or 
\begin{equation}\label{sec5-eq2}
(\psi(u_1),\psi(u_2),
\psi(u\idf{3}{4}),\psi(u_5)) 
=(v_3, v_2,v\idf{1}{5},v_4).
\end{equation}

Thus it is clear that 
the selection of $\phi$ and $\psi$ 
above may be not unique. 
%It is possible that $\phi$ is also an isomorphism from $G$ to $H$.

\item[Step 3.]
Let $D_{\psi}$ denote the 
digraph with vertex set $[n]$ such that $i\rightarrow j$ is an arc 
of this digraph 
if and only if one of the following cases occurs:
\begin{enumerate}
\item[(a)] $i\notin \{s_1,s_2\}, j\notin \{t_1,t_2\}$ and 
$\psi(u_i)=v_j$;
\item[(b)] $i\in \{s_1,s_2\}$, $j\notin \{t_1,t_2\}$ and $\psi(u\idf{s_1}{s_2})=v_j$;

\item[(c)] $i\notin \{s_1,s_2\}$, $j\in \{t_1,t_2\}$ and $\psi(u_i)=v\idf{t_1}{t_2}$; and

\item[(d)] $i\in \{s_1,s_2\}$, $j\in \{t_1,t_2\}$ and $\psi(u\idf{s_1}{s_2})=v\idf{t_1}{t_2}$.
\end{enumerate}
%Note that the condition ``$\psi$ is an isomorphism from $G/e$ to $H/f$" is not applied in the definition of $D_{\psi}$.
Also observe that the digraph $D_{\psi}$ 
has exactly $n+2$ arcs if
$\psi(u\idf{s_1}{s_2})=v\idf{t_1}{t_2}$,
%$\psi(\{u_{s_1},u_{s_2}\})=\{v_{t_1},v_{t_2}\}$,
and %has  exactly 
$n+1$ arcs otherwise. 

For example, if $n=5$, $s_1=3,s_2=4,t_1=1,t_2=5$ and 
$$
\psi(u_1)=v_4, \psi(u_2)=v\idf{1}{5},
%\{v_1,v_5\}, \psi(\{u_3,u_4\})
\psi(u\idf{3}{4})=v_2, \psi(u_5)=v_3,
$$
then the digraph $D_{\psi}$ is 
the one %the digraph  
shown in Figure~\ref{f-G10-5},
which contains two directed cycles 
$C_1: 1\rightarrow 4\rightarrow 2\rightarrow 1$
and $C_2: 2\rightarrow 5\rightarrow 3\rightarrow 2$.

%The digraph $D_{\psi}$  has the property in Proposition~\ref{dig1}.

\tikzstyle{cblue}=[circle, draw, thin,fill=blue!20, scale=0.5]
\begin{figure}[!ht]
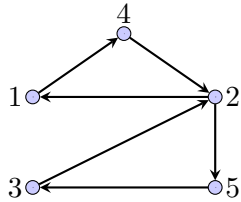

	\tikzp{1.2}
	{
		\foreach \place/\y in {{(0,0)/1}, {(2,0)/2},{(1,0.7)/3},{(2,-1)/4},{(0,-1)/5}}   
		\node[cblue] (b\y) at \place {};
		
		\filldraw[black] (b1) circle (0pt)node[anchor=east] {1};
		\filldraw[black] (b2) circle (0pt)node[anchor=west] {2};
		\filldraw[black] (b3) circle (0pt)node[anchor=south] {4};
		\filldraw[black] (b4) circle (0pt)node[anchor=west] {5};
		\filldraw[black] (b5) circle (0pt)node[anchor=east] {3};		

		\draw[black, ->, >=stealth, thick] (b2)--(b4);
		\draw[black, ->, >=stealth, thick] (b4)--(b5);
		\draw[black, ->, >=stealth, thick] (b5)--(b2);
		\draw[black, ->, >=stealth, thick] (b2)--(b1);				
		\draw[black, ->, >=stealth, thick] (b1)--(b3);			
		\draw[black, ->, >=stealth, thick] (b3)--(b2);	
											}

 \caption{\small Digraph $D_{\psi}$ with two directed cycles
$1\rightarrow 4\rightarrow 2\rightarrow 1$
and $2\rightarrow 5\rightarrow 3\rightarrow 2$}
\label{f-G10-5}
\end{figure}

\item[Step 4.] Let $C_1, C_2,\dots, C_r$ be the 
distinct directed cycles in $D_{\psi}$,
where $r\ge 1$ and $C_i$ denotes the cycle
$\pi_i(1)\rightarrow \pi_i(2)\rightarrow 
\cdots \rightarrow \pi_i(k_i)\rightarrow \pi_i(1)$
in $D_{\psi}$.

For each $i\in [r]$, %$i=1,2,\dots,r$, 
let $W_i$ be any graph with a vertex orbit 
$(w_{i,1},\dots,w_{i,k_i})$ of 
some automorphism $\xi_i$
such that $\xi_i(w_{i,j})=w_{i,j+1}$ for all $j\in [k_i]$,
where $w_{i,k_i+1}$ and $w_{i,1}$ are the same vertex
in $W_i$.

Let $G_0:=G$ and $H_0:=H$, 
and for each $i\in [r]$, 
let $G_i$ (resp. $H_i$) denote
the graph obtained from $G_{i-1}$ (resp. $H_{i-1}$)
and $W_i$ 
by identifying $u_{\pi_i(j)}$ (resp. $v_{\pi_i(j)}$) 
and $w_{i,j}$ 
for all $j\in [k_i]$.
It will be shown in Theorem~\ref{new-member} that 
$(G_r,e,H_r,f)$ is a member of $\Phi'$.

\iffalse
Let $G'$ (resp. $H'$) denote
the graph obtained from $G$ (resp. $H$) 
and $W_1,W_2,\dots,W_r$ 
by identifying $u_{\pi_i(j)}$ (resp. $v_{\pi_i(j)}$) 
and $w_{i,j}$ 
for all $i=1,2,\dots,r$ and $j=1,2,\dots,k_i$.
\fi

For example, for graphs $G_0$ and $H_0$ in 
Figure~\ref{f-G10-2-0}, 
if $\psi$ is the mapping defined in  (\ref{sec5-eq1}),
then $D_{\psi}$ is the digraph shown in 
Figure~\ref{f-G10-5} which contains two directed cycles 
$C_1: 1\rightarrow 4\rightarrow 2\rightarrow 1$
and $C_2: 2\rightarrow 5\rightarrow 3\rightarrow 2$.

Let $W_1$ be the graph $K_3$ with vertices 
$w_{1,1},w_{1,2}$ and $w_{1,3}$,
and let $W_2$ be the graph $K_{1,3}$ 
with $w_{2,1},w_{2,2}$ and $w_{2,3}$ be its three leaves,
as shown below.

\tikzstyle{cblue}=[circle, draw, thin,fill=blue!20, scale=0.5]
\begin{figure}[!ht]
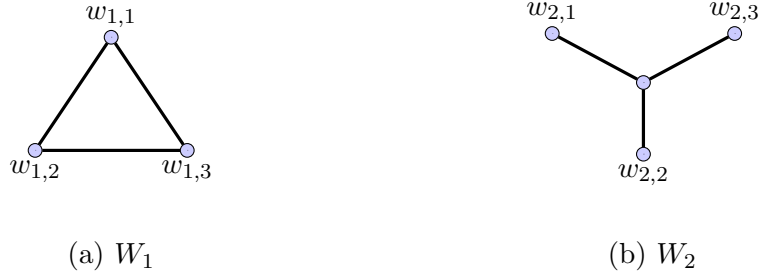

	\tikzp{1}
	{
		\foreach \place/\y in {{(1,-0.75)/1}, {(-1,-0.75)/2},{(0,0.75)/3}}   
		\node[cblue] (b\y) at \place {};
		
		\filldraw[black] (b1) circle (0pt)node[anchor=north] {$w_{1,3}$};
		\filldraw[black] (b2) circle (0pt)node[anchor=north] {$w_{1,2}$};
		\filldraw[black] (b3) circle (0pt)node[anchor=south] {$w_{1,1}$};

		\draw[black, very thick] (b1) -- (b2)--(b3)--(b1);

		\foreach \place/\y in {{(0+7,0.15)/1}, {(-1.2+7,0.8)/2},{(1.2+7,0.8)/3}, {(0+7,-0.8)/4}}   
		\node[cblue] (c\y) at \place {};
		
		\filldraw[black] (c1) circle (0pt)node[anchor=south] {};
		\filldraw[black] (c2) circle (0pt)node[anchor=south] {$w_{2,1}$};
		\filldraw[black] (c3) circle (0pt)node[anchor=south] {$w_{2,3}$};
		\filldraw[black] (c4) circle (0pt)node[anchor=north] {$w_{2,2}$};

		\draw[black, very thick] (c3)-- (c1) -- (c2);	 
		\draw[black, very thick]  (c1) -- (c4);	 						}

	{}\hfill \hspace{-0.2cm} (a) $W_1$ \hspace{5.7 cm} (b) $W_2$\hfill {}			
 \caption{\small Graphs $W_1$ and $W_2$}
\label{f-G10-6}
\end{figure}

Clearly, for $i=1,2$, $W_i$ has an 
automorphism $\xi_i$
such that $\xi_i(w_{i,j})=w_{i,j+1}$ for all $j=1,2,3$,
where $w_{i,4}$ is the vertex $w_{i,1}$ for $i=1,2$.

Now, let $G_2$ be the graph obtained from $G_0$ 
and $W_1, W_2$ by identifying each pair of the following 
vertices:

$u_1$ and $w_{1,1}$; $u_4$ and $w_{1,2}$; $u_2$ and $w_{1,3}$;

$u_2$ and $w_{2,1}$; $u_5$ and $w_{2,2}$; $u_3$ and $w_{2,3}$.

Similarly, let $H_2$ be the graph obtained from $H_0$ 
and $W_1, W_2$ by identifying each pair of the following 
vertices:

$v_1$ and $w_{1,1}$; $v_4$ and $w_{1,2}$; $v_2$ and $w_{1,3}$;

$v_2$ and $w_{2,1}$; $v_5$ and $w_{2,2}$; $v_3$ and $w_{2,3}$.

Then $G_2$ and $H_2$ are the graphs shown below, which  are exactly the examples given in Figure~\ref{f-G10}.

\tikzstyle{cblue}=[circle, draw, thin,fill=blue!20, scale=0.5]
\begin{figure}[!ht]
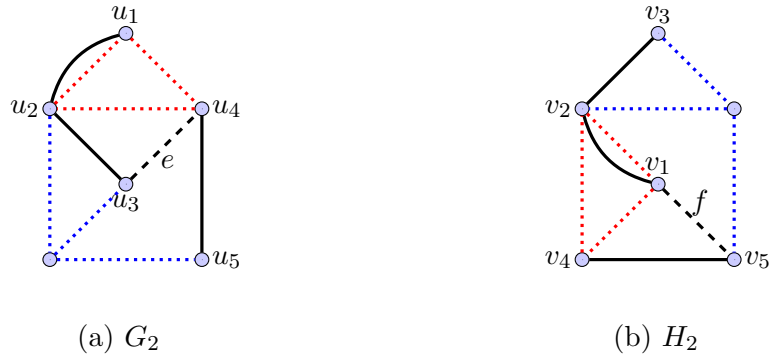

	\tikzp{1}
	{
		\foreach \place/\y in {{(0,0)/1}, {(1,1)/2},{(1,-1)/3}, {(-1,-1)/4}, {(-1,1)/5}, {(0,2)/6}}   
		\node[cblue] (b\y) at \place {};
		
		\filldraw[black] (b1) circle (0pt)node[anchor=north] {$u_3$};
		\filldraw[black] (b2) circle (0pt)node[anchor=west] {$u_4$};
		\filldraw[black] (b3) circle (0pt)node[anchor=west] {$u_5$};
		\filldraw[black] (b4) circle (0pt)node[anchor=south] {};
		\filldraw[black] (b5) circle (0pt)node[anchor=east] {$u_2$};
		\filldraw[black] (b6) circle (0pt)node[anchor=south] {$u_1$};

				\node [style=none] (cap1) at (0.55,0.3)  {$e$};
				
		\draw[black, very thick] (b5) -- (b1);
		\draw[black, very thick] (b2) -- (b3);		
		\draw[black, red, dotted, very thick] (b2) -- (b5) -- (b6) -- (b2);		
		\draw[black, blue, dotted, very thick] (b5) -- (b4) -- (b1);		
		\draw[black, blue, dotted, very thick] (b4) -- (b3);		
		\draw[black, dashed, very thick] (b1) -- (b2);								
		\draw [black, very thick, bend left=30] (b5) to (b6);

		\foreach \place/\y in {{(0+7,0)/1}, {(1+7,1)/2},{(1+7,-1)/3}, {(-1+7,-1)/4}, {(-1+7,1)/5}, {(0+7,2)/6}}   
		\node[cblue] (c\y) at \place {};
		
		\filldraw[black] (c1) circle (0pt)node[anchor=south] {$v_1$};
		\filldraw[black] (c2) circle (0pt)node[anchor=south] {};
		\filldraw[black] (c3) circle (0pt)node[anchor=west] {$v_5$};
		\filldraw[black] (c4) circle (0pt)node[anchor=east] {$v_4$};
		\filldraw[black] (c5) circle (0pt)node[anchor=east] {$v_2$};
		\filldraw[black] (c6) circle (0pt)node[anchor=south] {$v_3$};

				\node [style=none] (cap1) at (0.55+7,-0.25)  {$f$};
				
		\draw[black, very thick] (c4) -- (c3);	 
		\draw[black, very thick] (c5) -- (c6);			
		\draw[black, very thick, dashed] (c1) -- (c3);
		\draw [black, very thick, bend right=30] (c5) to (c1);
		\draw[black, red, dotted, very thick] (c5) -- (c1) -- (c4) -- (c5);		
		\draw[black, blue, dotted, very thick] (c6) -- (c2) -- (c5);		
		\draw[black, blue, dotted, very thick] (c2) -- (c3);				
						}

	{}\hfill \hspace{-0.2cm} (a) $G_2$ \hspace{5.7 cm} (b) $H_2$\hfill {}			
 \caption{\small Graphs $G_2$ and $H_2$}
\label{f-G10-7}
\end{figure}

\iffalse 
\begin{figure}[htbp]
	\centering
	
	\input f-G10-7-1.pic
	
	\vspace{0.3 cm}
	
	\hspace{0.2cm} $G_2$ \hspace{4.8 cm} $H_2$
	
	\caption{\small Graphs $G_2$ and $H_2$}
	\label{f-G10-7-1}
\end{figure}
\fi 

\iffalse 
For any directed cycle 
$\pi(1)\rightarrow \pi(2)\rightarrow 
\cdots \rightarrow \pi(k)\rightarrow \pi(1)$ 
in $D_{\psi}$, where $k\ge 1$, and 
any graph $W$ with an vertex orbit
$w_1,w_2,\dots,w_k$ of some automorphism $\xi$
such that $\xi(w_i)=w_{i+1}$ for all $i=1,2,\dots,k$,
where $w_{k+1}$ is assumed to be $w_1$, 
$(G',e,H',f)$ must be a member of $\Phi'$,
where  $G'$ and $H'$ denote graphs 
$G(u_{\pi(1)}, u_{\pi(2)},\dots,u_{\pi(k)})\sqcup 
W(w_1,w_2,\dots,w_k)$ 
and $H(v_{\pi(1)}, v_{\pi(2)},\dots,v_{\pi(k)})\sqcup 
W(w_1,w_2,\dots,w_k)$ respectively.
\fi 
%%%%%%%%%

%\item $(G',e,H,f)\in \Phi'$.
\end{enumerate}

Now we prove that the above construction fulfills our goal based on the following observation.

\begin{pro}\label{dig1}
Assume that $C$ is any directed cycle in the digraph
$D_{\psi}$. 
Then both $s_1$ and $s_2$ are on $C$  
if and only if both $t_1$ and $t_2$ are on $C$.
If both $s_1$ and $s_2$ are on $C$,
then   
$\psi(u\idf{s_1}{s_2})
=v\idf{t_1}{t_2}$
%$\psi(\{u_{s_1},u_{s_2}\})=\{v_{t_1},v_{t_2}\}$
and $C$  contains both arcs $s_1\rightarrow t_j$ and 
$s_2\rightarrow t_{3-j}$ for some $j\in \{1,2\}$.
\end{pro}

\proof 
Assume that both $s_1$ and $s_2$ are on $C$. 
If $\psi(u\idf{s_1}{s_2})=v_q$
for some $q\in [n]$, 
then, $s_1\rightarrow q$ 
(resp. $s_2\rightarrow q$) 
is the only arc in $D_{\psi}$
leaving $s_1$ (resp. $s_2$),
implying that $C$ contains both 
arcs $s_1\rightarrow q$ and $s_2\rightarrow q$, a 
contradiction to the fact that $C$
is a directed cycle in $D_{\psi}$.
\iffalse 
As $C$ is a directed cycle, 
for any $q\in [n]$,
$C$ cannot  contain both arcs 
$s_1\rightarrow q$ and $s_2\rightarrow q$,
%for any $q\in [n]$,
 implying that 
$\psi(u\idf{s_1}{s_2})
%\{u_{s_1},u_{s_2}\})
\ne v_q$.\fi 
Thus
$\psi(u\idf{s_1}{s_2})
=v\idf{t_1}{t_2}$
%$\psi(\{u_{s_1},u_{s_2}\})=\{v_{t_1},v_{t_2}\}$
by the definition of $\psi$.
As $s_1,s_2$ are on $C$ and 
$C$ cannot contain both arcs 
$s_1\rightarrow t_j$ and $s_2\rightarrow t_j$
for any $j\in \{1,2\}$, 
$C$ must contain both arcs $s_1\rightarrow t_j$ and 
$s_2\rightarrow t_{3-j}$ for some $j\in \{1,2\}$.
Hence $t_1,t_2$ are also contained on $C$.

Similarly, it can be shown that 
if both $t_1$ and $t_2$ are on 
a directed cycle $C$,
then
$\psi(u\idf{s_1}{s_2})=v\idf{t_1}{t_2}$, 
implying that 
$s_1,s_2$ are also on $C$ 
which contains both arcs $s_1\rightarrow t_j$ and 
$s_2\rightarrow t_{3-j}$ for some $j \in [2]$.
\proofend

\begin{theo}\label{new-member}
For the graphs $G_r$ and $H_r$ obtained in Step 4 above, 
we have  $(G_r,e,H_r,f)\in \Phi'$.
\iffalse 
Let 
$\pi(1)\rightarrow \pi(2)\rightarrow 
\cdots \rightarrow \pi(k)\rightarrow \pi(1)$ 
be any directed cycle in $D_{\psi}$, where $k\ge 1$, and 
$W$ be any graph with an vertex orbit
$\{w_i: i\in [k]\}$ of some automorphism $\xi$,
where $\xi(w_i)=w_{i+1}$ for all $i=1,2,\dots,k$
and $w_{k+1}$ and $w_1$ are assumed to be the same vertex in $W$.
Then  $(G',e,H',f)\in \Phi'$,
where  $G'$ and $H'$ denote graphs 
$G(u_{\pi(1)}, u_{\pi(2)},\dots,u_{\pi(k)})\sqcup 
W(w_1,w_2,\dots,w_k)$ 
and $H(v_{\pi(1)}, v_{\pi(2)},\dots,v_{\pi(k)})\sqcup 
W(w_1,w_2,\dots,w_k)$ respectively.
\fi
\end{theo}

\proof We shall show that $(G_i,e,H_i,f)\in \Phi'$ for each $i\in[r]$.
Let $G_0$ and $H_0$ be the graphs
$G$ and $H$ respectively, and 
let $\phi_0$ be $\phi$ and $\psi_0$ 
be $\psi$. 
It suffices to prove the two results below:

\begin{enumerate}
	\item[(a)] For each $i\in [r]$, 
	if $(G_{i-1},e,H_{i-1},f)\in \Phi'$
	with an isomorphism $\phi_{i-1}$ from $G_{i-1}\backslash e$ to $H_{i-1}\backslash f$ and an isomorphism $\psi_{i-1}$ from 
	$G_{i-1}/ e$ to $H_{i-1}/ f$,
	then 
	$(G_i,e,H_i,f)\in \Phi'$
	with an isomorphism $\phi_i$ from $G_i\backslash e$ to $H_i\backslash f$ and an isomorphism $\psi_i$ from $G_i/ e$ to $H_i/ f$,
	where 
	$\phi_i$ is  extended from $\phi_{i-1}$ by defining $\phi_i(w)=w$ for each  $w\in V(W_i)\setminus \{w_{i,j}: j\in [k_i]\}$ and
	$\psi_i$ is extended from $\psi_{i-1}$ by defining $\psi_i(w)=\xi_i(w)$ for each $w\in V(W_i)\setminus \{w_{i,j}: j\in [k_i]\}$.
	
\item[(b)] For any $2\le i\le r$, 
$C_i$ is a directed cycle in the digraph $D_{\psi_{i-1}}$.
\end{enumerate} 
Notice that  results (a) and (b) above
can be obtained by applying 
the following conclusions (i) and (ii)
repeatedly: 
\begin{enumerate}
\item $(G_1,e,H_1,f)\in \Phi'$, as 
 $\phi_1$ is an isomorphism from $G_1\backslash e$ to $H_1\backslash f$ and  $\psi_1$ 
 is an isomorphism  from $G_1/ e$ to $H_1/ f$, and 
\iffalse 
where 
$\phi_1$ is  extended from $\phi$ by defining $\phi_1(w)=w$ for each  $w\in V(W_1)\setminus \{w_{1,1},\dots,w_{1,k_1}\}$ and
$\psi_1$ is extended from $\psi$ by defining $\psi_1(w)=\xi_1(w)$ for each $w\in V(W_1)\setminus \{w_{1,1},\dots,w_{1,k_1}\}$.\fi

\iffalse 
 \red{where $\phi_1$ is the isomorphism from $G_1\backslash e$ to $H_1\backslash f$ extended from $\phi$ by defining $\phi_1(w)=w$ for each  $w\in V(W_1)\setminus \{w_{1,1},\dots,w_{1,k_1}\}$ and
 $\psi_1$ is the isomorphism from $G_1/e$ to $H_1/f$ extended from $\psi$ by defining $\psi_1(w)=\xi_1(w)$ for each $w\in V(W_1)\setminus \{w_{1,1},\dots,w_{1,k_1}\}$.}
\fi 
 %from $(V(G_1)-\{u_{s_1},u_{s_2}\})\cup \{\{u_{s_1},u_{s_2}\}\}$ to $(V(H_1)-\{v_{t_1},v_{t_2}\})\cup \{\{v_{t_1},v_{t_2}\}\}$ extended from $\psi$ by defining $\psi_1(w)=\xi_1(w)$ for each $w\in V(W_1)\setminus \{w_{1,1},\dots,w_{1,k_1}\}$.
%the theorem holds for $r=1$;

\item if $r\ge 2$, then $C_2$ is a directed cycle in the digraph $D_{\psi_1}$.
%, where $\psi_1$ is the bijection from $(V(G_1)-\{u_{s_1},u_{s_2}\})\cup \{\{u_{s_1},u_{s_2}\}\}$ to $(V(H_1)-\{v_{t_1},v_{t_2}\})\cup \{\{v_{t_1},v_{t_2}\}\}$ extended from $\psi$ by defining $\psi_1(w)=\xi_1(w)$ for each $w\in V(W_1)\setminus \{w_{1,1},\dots,w_{1,k_1}\}$.

%if $G'$ and $H'$ are the graphs obtained in Step 4
%with $r=1$ and the directed cycle $C_1$ in $D_{\psi}$,  
%then any directed cycle $C_2$ in $D_{\psi}$ distinct from $C_1$ 
%is also a directed cycle in the digraphs $D_{\psi'}$,
%where $\psi'$ is the bijection from 
%$(V(G')-\{u_{s_1},u_{s_2}\})\cup \{\{u_{s_1},u_{s_2}\}\}$ to 
%$(V(H')-\{v_{t_1},v_{t_2}\})\cup \{\{v_{t_1},v_{t_2}\}\}$ 
%extended from $\psi$ by defining 
%$\psi'(w)=\xi_1(w)$ for each $w\in V(G_1)$.
\end{enumerate}

For (i), it is easy to verify that $\phi_1$ is an isomorphism from $G_1\backslash e$ to $H_1\backslash f$. In the following, we shall prove that  $\psi_1$ is an isomorphism from $G_1/e$ to $H_1/f$, i.e., $\psi_1$ is a bijection from $V(G_1/e)$ to $V(H_1/f)$ such that 
%$xy\in E(G_1/e)$ if and only if $\psi_1(x)\psi_1(y)\in E(H_1/f)$.
$\epsilon_{G_1/e}(x,y)=\epsilon_{H_1/f}(\psi_1(x),\psi_1(y))$ holds for each pair of vertices $x,y\in V(G_1/e)$,
where $\epsilon_G(u,v)$ is the number of edges in $G$ 
joining $u$ and $v$.
Since $\psi_1|_{V(G/e)}=\psi$ is an isomorphism from $G/e$ to $H/f$, this conclusion is obvious
if both $x$ and $y$ are contained in 
$V(G/e)$.
Thus, we assume that 
$y\notin
V(G/e)$.

\noindent {\bf Case 1}: 
$\{s_1,s_2\}\not\subseteq V(C_1)$.

By Proposition~\ref{dig1}, $\{t_1,t_2\}\not\subseteq V(C_1)$ in this case. Thus, 
$$
V(G_1/e)\setminus (V(G/e)\setminus \{u_{\pi_1(j)}: j\in [k_1]\})
%,\dots,u_{\pi_1(k_1)}\})
=V(W_1)
=V(H_1/f)\setminus (V(H/f)\setminus
\{v_{\pi_1(j)}: j\in [k_1]\}) %\{v_{\pi_1(1)},\dots,v_{\pi_1(k_1)}\}).
$$
Due to the construction in Step 4, for any $w_{1,s}\in V(W_1)$, 
$$
\psi_1(w_{1,s})=\psi(u_{\pi_1(s)})=v_{\pi_1(s+1)}=w_{1,s+1}=\xi_1(w_{1,s}),
$$ 
implying that $\psi_1|_{V(W_1)}=\xi_1$. Thus $\psi_1$ is a bijection as $\xi_1$ is an automorphism of $W_1$.
Also, for any $x\in V(G_1/e)$ and $y\in V(G_1/e)\setminus V(G/e)$, $xy\in E(G_1/e)$ only when $x\in V(W_1)$, and  
$\epsilon_{G_1/e}(x,y)=\epsilon_{W_1}(\xi_1(x),\xi_1(y))=\epsilon_{H_1/f}(\psi_1(x),\psi_1(y))$.
%$xy\in E(G_1/e)$  if and only if $\psi_1(x)\psi_1(y)=\xi_1(x)\xi_1(y)\in E(H_1/f)$. 
Hence $\psi_1$ is an isomorphism in this case.
%. Note that in this case, $xy\in E(G_1/e)$ only when $x\in \left(V(G_1/e)\setminus V(G/e)\right)\cup \{w_{1,1},\dots,w_{1,k_1}\}$. Due to the construction in Step 4 and the automorphism $\xi_1$ of $W_1$, $\psi_1$ is a bijection, and $xy\in E(G_1/e)$ if and only if $\psi_1(x)\psi_1(y)=\xi_1(x)\xi_1(y)\in E(H_1/f)$.

\noindent {\bf Case 2}: 
$\{s_1,s_2\}\subseteq V(C_1)$.

%Otherwise, $\{s_1,s_2\}\subseteq C_1$, and 
By Proposition~\ref{dig1}, $\{t_1,t_2\}\subseteq V(C_1)$.
Assume that $w_{1,p}$ and $w_{1,q}$ are identified with $u_{s_1}$ and $u_{s_2}$ in the construction of $G_1$, respectively. Then, by Proposition~\ref{dig1}, $\xi_1(w_{1,p})$ and $\xi_1(w_{1,q})$ are identified with vertices $v_{t_j}$ and $v_{t_{3-j}}$ in the construction of $H_1$ respectively, for some $j\in [2]$. 
As a result, 
$$
V(G_1/e)\setminus (V(G/e)\setminus
(\{u_{\pi_1(j)}:j\in [k_1]\}\cup 
\{u\idf{s_1}{s_2}\}) )
 %\{u_{\pi_1(1)},\dots,u_{\pi_1(k_1)},u\idf{s_1}{s_2} \})
=V(W_1\cdot w_{1,p}w_{1,q})
$$ 
and  
$$
V(H_1/f)\setminus (V(H/f)\setminus
(\{v_{\pi_1(j)}:j\in [k_1]\}\cup 
\{v\idf{t_1}{t_2}\}) )
% \{v_{\pi_1(1)},\dots,v_{\pi_1(k_1)},v\idf{t_1}{t_2}%\{v_{t_1},v_{t_2}\}\})
=V(W_1\cdot \xi_1(w_{1,p})\xi_1(w_{1,q})).
$$

It is clear to see that $\xi_1$ induces an isomorphism from $W_1\cdot w_{1,p}w_{1,q}$ to $W_1\cdot \xi_1(w_{1,p})\xi_1(w_{1,q})$ naturally, which implies that $\psi_1$ is a bijection, and for any $x\in V(G_1/e)$ and $y\in V(G_1/e)\setminus V(G/e)$,
%$xy\in E(G_1/e)$ if and only if $\psi_1(x)\psi_1(y)=\xi_1(x)\xi_1(y)\in E(H_1/f)$ 
$$\epsilon_{G_1/e}(x,y)=\epsilon_{W_1\cdot w_{1,p}w_{1,q}}(x,y)=\epsilon_{W_1\cdot \xi_1(w_{1,p})\xi_1(w_{1,q})}(\xi_1(x),\xi_1(y))=\epsilon_{H_1/f}(\psi_1(x),\psi_1(y)).$$ Thus $\psi_1$ is an isomorphism. (i) holds.

\iffalse$\xi_1'$ be the mapping from $V(W_1\cdot w_{1,p}w_{1,q})$ to $V(W_1\cdot \xi_1(w_{1,p})\xi_1(w_{1,q}))$ by defining $\xi_1'(w)=\xi_1(w)$ for all $w\in W_1$. Then $\xi_1'$ is a well-defined isomorphism. 
In other words, $w_{1,p}$ is the same vertex as $w_{1,q}$ in graph $G_1/e$ and $\xi_1(w_{1,p})$ is the same vertex as $\xi_1(w_{1,q})$ in graph $H_1/f$.
\fi

Since $\psi_1|_{V(G/e)}=\psi$, (ii) holds.
The proof is complete.
\proofend

\noindent\textbf{Remark:} Similar non-isomorphic $T$-equivalent pairs can be obtained repeatedly through Steps 1--4 starting from any path longer than 3 with respect to its second rightmost edge.

\section*{Acknowledgement}

This research is supported 
by NSFC (No.  12371340)
and 
the Ministry of Education,
Singapore, under its Academic Research Tier 1 (RG19/22). Any opinions,
findings and conclusions or recommendations expressed in this
material are those of the author(s) and do not reflect the views of the
Ministry of Education, Singapore.

\end{document}